\documentclass[11pt,a4paper]{amsart}
\usepackage[T1]{fontenc}
\usepackage{newpxtext,newpxmath}    
\usepackage{graphicx}
\usepackage[english]{babel}
\usepackage{amsthm,mathtools,nicefrac}
\usepackage[margin=0.9in]{geometry}
\usepackage[shortlabels]{enumitem}  
\usepackage{etoolbox}
\usepackage{xcolor}
\usepackage{todonotes}
\usepackage{stmaryrd}
\DeclareFontFamily{U}{stmry}{}
\DeclareFontShape{U}{stmry}{m}{n}{<-> stmary10}{}
\usepackage{dsfont} 
\usepackage{tikz}
\usetikzlibrary{calc,arrows.meta}
\usepackage[colorlinks=true,linkcolor=blue!60!black,citecolor=green!40!black,urlcolor=blue!60!black]{hyperref}

\newcommand{\By}[2]{\overset{\mbox{\tiny{#1}}}{#2}}
\newcommand{\ByRef}[2]{   \By{\eqref{#1}}{#2} }
\newcommand{\eqBy}[1]{    \By{#1}{=} }
\newcommand{\lBy}[1]{     \By{#1}{<} }

\newcommand{\leBy}[1]{    \By{#1}{\le} }

\newcommand{\eqByRef}[1]{ \ByRef{#1}{=} }

\newcommand{\leByRef}[1]{ \ByRef{#1}{\le} }
\newcommand{\geByRef}[1]{ \ByRef{#1}{\ge} }

\DeclareMathOperator{\dist}{dist}
\newcommand{\nb}{w}

\theoremstyle{plain}
\newtheorem{lem}{Lemma}[section]
\newtheorem{fact}[lem]{Fact}

\newtheorem{defi}[lem]{Definition}
\newtheorem{rem}[lem]{Remark}
\newtheorem{thm}[lem]{Theorem}

\newtheorem{prop}[lem]{Proposition}

\newcounter{claim}
\renewcommand{\theclaim}{\Alph{claim}}

\newcounter{proofctr}

\pretocmd{\proof}{%
  \stepcounter{proofctr}%
  \setcounter{claim}{0}%
}{}{}

\makeatletter
\newcommand{\labelleditem}[2]{\item[\textup{#1}]\phantomsection\protected@edef\@currentlabel{\textup{#1}}\label{#2}}
\makeatother

\newcommand{\eps}{\varepsilon}

\DeclarePairedDelimiter{\norm}{\lVert}{\rVert}

\makeatletter
\patchcmd{\abstract}{\item[\hskip\labelsep\scshape\abstractname.]}%
  {\item[]\makebox[\linewidth][c]{\scshape\abstractname}\par\nobreak\vspace{10pt}\noindent\ignorespaces}%
  {}{\PackageWarning{etoolbox}{Patching the abstract heading failed}}
\patchcmd{\@setabstracta}{\skip@20\p@}{\skip@28\p@}%
  {}{\PackageWarning{etoolbox}{Patching the abstract spacing failed}}
\newcommand{\affiliation}[1]{\gdef\@affiliation{#1}}
\apptocmd{\@setauthors}{%
  \ifdefined\@affiliation
    {\par\vspace{2pt}\centering\footnotesize\@affiliation\par}%
  \fi}{}{}
\makeatother

\author[F.~Garbe]{Frederik Garbe}
\affiliation{Institute of Computer Science\\ Czech Academy of Sciences}
\thanks{Research supported by Czech Science Foundation Project 26-23695S}

\title[An improved Moore bound]{An improved Moore bound for girths seven and eleven}
\address{Institute of Computer Science of the Czech Academy of Sciences, Prague, Czech Republic}
\email{garbe@cs.cas.cz}

\begin{document}

\begin{abstract}
Let $n(d,k)$ be the minimum number of vertices of a graph of minimum degree $d$ and girth at least $g=2k+1$. A simple counting argument gives the classical lower bound $n(d,k)\ge m(d,k)=1+d\sum_{i=0}^{k-1}(d-1)^i$, known as the \emph{Moore bound} and dating back at least to 1960. Apart from complete graphs and odd cycles, graphs attaining this lower bound can exist only for girth $5$ and degrees $3,7$ and $57$ (Hoffman and Singleton 1960; Bannai and Ito 1973; Damerell 1973). Nevertheless, for more than sixty years this lower bound has not been substantially improved beyond an additive constant. For $k\geq 3$, i.e. girth at least $7$, the best general bound was $n(d,k)\ge m(d,k)+2$, due to Bannai and Ito (1981). We show that for fixed girth $7$ and growing degree $d$ we have $n(d,3)\ge m(d,3)+(1-o(1))d/8$. Furthermore, for girth $11$ we have $n(d,5)\ge m(d,5)+\Omega(d^{3/5})$. In particular, the \emph{excess} $n(d,k)-m(d,k)$, for $k\in\{3,5\}$, tends to infinity with $d$ and is therefore generally unbounded, which answers a question of Bollob\'as from 1978.
For the degree/diameter problem we obtain similar bounds. A graph of maximum degree $d$ and diameter $3$ or $5$ has at most $m(d,3)-\Omega(\sqrt{d})$ or $m(d,5)-\Omega(d^{3/10})$ vertices, respectively. This answers a question of Bermond and Bollob\'as from 1981: for every constant $c$ there are $d$ and $k$ such that every graph of maximum degree $d$ and diameter $k$ has at least $c$ fewer vertices than the Moore bound.
\end{abstract}

\maketitle

\section{Introduction}

Let $G$ be a graph of girth at least $g=2k+1$ and minimum degree $d$. Starting a breadth-first search at an arbitrary vertex shows that $G$ contains as a subgraph the tree of depth $k$ in which every non-leaf vertex has degree $d$. This gives the lower bound $m(d,k)\coloneqq 1+d\sum_{i=0}^{k-1}(d-1)^i$ on the number of vertices of $G$, which is called the \emph{Moore bound}. Alon, Hoory and Linial~\cite{AlHoLi02} showed that this bound remains valid even if the minimum degree is replaced by the average degree, answering an old problem of Bollob\'as~\cite[p.~163, Problem~10]{Bol78}.
Graphs which attain this bound are always regular and are called \emph{Moore graphs}.
Trivial examples of such graphs are complete graphs (girth $3$) and odd cycles (degree $2$).
Hoffman and Singleton~\cite{HoSi60} showed that there are no such graphs of girth $7$ and $d\ge 3$, and that for girth $5$ Moore graphs have degree $2,3,7$ or $57$. The existence of a Moore graph of degree $57$ is still a well-known open problem. Later, Bannai and Ito~\cite{BaIt73} and independently Damerell~\cite{Dam73} showed that there are no Moore graphs of degree $d\ge 3$ and odd girth at least $9$. The same question can be investigated for even girth $g=2k$ and leads to an analogous bound, which for $d\ge3$ can only be attained for girths $4,6,8$ and $12$ by results of Feit and Higman~\cite{FeHi64} and Singleton~\cite{Sin66}. For each of these girths, incidence graphs of generalised polygons attain it whenever $d-1$ is a prime power~\cite{Ben66,VMa98}. However, in the following we will restrict the discussion to the odd case and always assume that $g=2k+1$ unless explicitly said otherwise.

It is a natural extremal question to ask for the minimum number of vertices a graph with girth at least $g$ and minimum degree $d$ can have. We denote this number by $n(d,k)$.
For a graph $G$ of girth at least $g$ and minimum degree $d$ the difference $|V(G)|-m(d,k)$ is called the \emph{excess} $e$ of $G$. Note that if $e<(d-1)^k$ then the girth is exactly $2k+1$ and if $e<\sum_{i=0}^{k-1}(d-1)^i$ then $G$ must be a regular graph. Both are known facts which easily follow from the Moore bound. Our results only concern values in those ranges.

Brown~\cite{Bro67} and Bannai and Ito~\cite{BaIt81} showed that graphs with excess $1$ do not exist for any $k\ge 2$ and $d\ge 3$. Furthermore, Kov\'acs~\cite{Kov81} showed that graphs with excess $2$ and girth $5$ do not exist for a certain range of degrees. Lastly, the author showed in~\cite{Gar15} that graphs with girth $9$, excess $2$ and odd degree do not exist. This was the state of the art of the lower bound. In particular this means that the best known lower bound for the case $k=3$ was $n(d,3)\ge m(d,3)+2$ (ignoring computational results for very small values of $d$). In view of these improvements by only a small additive constant the question posed by Bollob\'as already in 1978~\cite[p.~163, Problem~11]{Bol78} feels even more natural: can the difference $n(d,k)-m(d,k)$ be bounded by an absolute constant?
Until now, the excess was known to be unbounded only in the special case of vertex-transitive graphs, by a result of Biggs~\cite{Big82}, and an answer to the original question for general graphs remained elusive.
We answer Bollob\'as's question negatively by achieving a substantial improvement of the Moore bound for girth $7$.

\begin{thm}\label{thm:girth7}
For $d\ge 3$, every graph of minimum degree $d$ and girth at least $7$ has more than $m(d,3)+(\sqrt d-2)^2/8$ vertices. In particular, a graph of minimum degree $d$, girth at least $7$ and excess $e$ satisfies $d<(\sqrt{8e}+2)^2$.
\end{thm}

Additionally, we show a similar result for girth $11$.

\begin{thm}\label{thm:girth11}
There exists $c>0$ such that for every $d\ge 3$ every graph of minimum degree $d$ and girth at least $11$ has more than $m(d,5)+c\,d^{3/5}$ vertices.
\end{thm}

Our proof uses eigenvalue considerations for the adjacency matrix of a graph near the Moore bound. This general approach dates back to the original paper of Hoffman and Singleton~\cite{HoSi60} and has since been used in different variations by most subsequent papers investigating the lower bound~\cite{Bro67,BaIt73,Dam73,BaIt81,Kov81,Gar15}. However, for more than sixty years it gave only improvements by an additive constant, and it seemed limited to these. Our paper is the first demonstration that it can give an unbounded improvement, in fact a linear one for girth~$7$. The crucial new idea, considering the shifted matrix $2A+I$ instead of $A$, is explained in the proof sketch below.

\subsection*{Degree/diameter problem}
There is a dual problem for which the Moore bound also gives the natural first bound. We denote by $\tilde{n}(d,k)$ the maximum number of vertices a graph with diameter at most $k$ and maximum degree $d$ can have. The same breadth-first search from above shows $\tilde{n}(d,k)\le m(d,k)$, as it reaches every vertex within $k$ steps and finds at most $d(d-1)^{i-1}$ new vertices in step $i$. For such a graph $G$, the difference $m(d,k)-|V(G)|$ is called the \emph{defect} of $G$. As for the degree/girth problem, only small constant improvements of this bound were known~\cite{FeMiPi11,Gar15}; see also the survey~\cite{MiSi13}. Again, this led to the question, this time of Bermond and Bollob\'as~\cite[Section~1]{BeBo81}, whether it is true that for every positive integer $c$ there are $k$ and $d$ such that $\tilde{n}(d,k)<m(d,k)-c$. In a complementary direction, Cames van Batenburg and Korsky~\cite{CvBK26} recently showed that $\tilde n(d,k)=(1+o(1))\,d^k$ for every fixed $k$, proving another conjecture of Bollob\'as. (For $k\in\{2,3,5\}$ this was known for infinitely many $d$~\cite{Del85a,Del85b}; see also~\cite{BSS19}.) So the smallest possible defect $m(d,k)-\tilde n(d,k)$ is known to be negligible compared with $m(d,k)$, and we show by the following theorems that it nevertheless tends to infinity for $k\in\{3,5\}$, answering the question of Bermond and Bollob\'as affirmatively. Previously, only conditional results were obtained for this problem~\cite{FiJa21}.

\begin{thm}\label{thm:diam3}
For $d\ge 3$, every graph of maximum degree $d$ and diameter at most $3$ has fewer than $m(d,3)-(\sqrt d-2)/\sqrt{8}$ vertices. In particular, a graph of maximum degree $d$, diameter at most $3$ and defect $f$ satisfies $d<(\sqrt{8}f+2)^2$.
\end{thm}

Similarly, we also show a result for diameter $5$.

\begin{thm}\label{thm:diam5}
There exists $c>0$ such that for every $d\ge 3$ every graph of maximum degree $d$ and diameter at most $5$ has fewer than $m(d,5)-c\,d^{3/10}$ vertices.
\end{thm}

\subsection*{Proof of a weaker result}
\begin{figure}[ht]
\centering
\begin{tikzpicture}[x=1cm,y=1cm,
  box/.style={draw=black!55, rounded corners=3pt, align=center, font=\small, inner sep=4pt, text width=3.55cm},
  dot/.style={circle, fill=black!60, inner sep=1.3pt}]
\draw[-{Stealth}, black!60] (-0.3,0) -- (14.2,0);
\node[box] (L) at (1.9,1.75) {\textbf{left outer cluster}\\ near $-\sqrt{2d}$\\ share $\approx\frac{1}{4}$};
\node[box, very thick] (M) at (6.6,1.75) {\textbf{middle cluster}\\ within $\frac{e}{d}$ of $-\frac{1}{2}$\\ share $\approx\frac{1}{2}$};
\node[box] (Rt) at (11.3,1.75) {\textbf{right outer cluster}\\ near $\sqrt{2d}$\\ share $\approx\frac{1}{4}$};
\foreach \x/\nm in {1.9/L, 6.6/M, 11.3/Rt} {\draw[black!55] (\nm.south) -- (\x,0) node[dot]{};}
\node[below=3pt, font=\small] at (1.9,0) {$-\sqrt{2d}$};
\node[below=3pt, font=\small] at (6.6,0) {$-\frac{1}{2}$};
\node[below=3pt, font=\small] at (11.3,0) {$\sqrt{2d}$};
\node[dot] at (13.5,0) {};
\node[below=3pt, font=\small] at (13.5,0) {$d$};
\node[above=3pt, font=\small, black!70] at (13.5,0) {once};
\end{tikzpicture}
\caption{The spectrum of a graph of girth at least $7$ and small excess.}
\label{fig:clusters}
\end{figure}
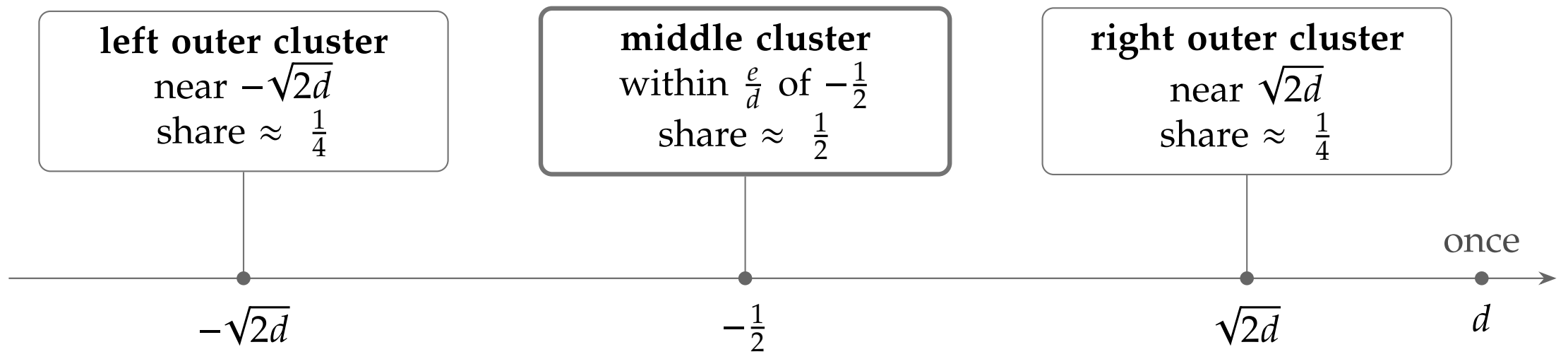

To demonstrate our method we prove a weaker version of Theorem~\ref{thm:girth7}, i.e. the bound
\begin{equation}\label{eq:weakerbound}
n(d,3)-m(d,3)=\Omega(\sqrt{d})\;.
\end{equation}
It is a short and elementary proof and nicely illustrates the key ideas of our paper while already giving a notable improvement to the Moore bound and answering the question of Bollob\'as. In the following we always assume that $d$ is sufficiently large.
Now, let $G$ be a graph with girth at least $7$ and minimum degree $d$ and excess
\begin{equation}\label{eq:exsqrtd}
e\le\sqrt{d}/7\;.
\end{equation}
Note that this implies that $G$ is connected, $d$-regular, and has girth exactly $7$. We define the \emph{excess graph} $E$ of $G$ to be the graph on the same vertex set in which two vertices are adjacent if they have distance at least $4$ in $G$. Note that $E$ is $e$-regular.
Let $A$ and $B$ be the adjacency matrices of $G$ and $E$, respectively. Furthermore, we define
\[F_d(x)\coloneqq x^3+x^2-2(d-1)x-(d-1)\;.\]
Recalling the well-known fact that the entry $(A^t)_{u,v}$ counts the number of walks of length $t$ from $u$ to $v$ in $G$ it is easy to see (using the fact that $G$ has girth $7$ and is $d$-regular) that
\begin{equation}\label{eq:matrixpoly}
F_d(A)=A^3+A^2-2(d-1)A-(d-1)I=J-B\;,
\end{equation}
where $I$ is the identity and $J$ is the all-one matrix. Since $G$ is $d$-regular and connected, $A$ has the eigenvalue $d$ with multiplicity $1$ and the eigenvector $\mathbf 1$. The eigenspace of every other eigenvalue $\lambda$ is orthogonal to $\mathbf 1$ and $\lambda\in[-d,d)$. Multiplying \eqref{eq:matrixpoly} from the right by an eigenvector corresponding to $\lambda$ gives $F_d(\lambda)=\mu_\lambda$ for an eigenvalue $-\mu_\lambda$ of $B$, for which also $|\mu_\lambda|\le e$. We call $\mu_\lambda$ the \emph{error} of $\lambda$.
So far, these are known observations about graphs with small excess, as previously used in~\cite{HoSi60,Bro67,Kov81,Gar15}. Now, to start our line of arguments we consider the distribution of those eigenvalues.
The polynomial $F_d(x)-\mu$, for $|\mu|\le e<d-1$, has degree $3$ and elementary calculations show that each of the following intervals contains one of the three roots (as the sign of the polynomial changes in each of those intervals):
\begin{equation}\label{eq:intervals}
I_-\coloneqq\left(-\sqrt{2(d-1)}-1,-\sqrt{2(d-1)}\right);\;I_0\coloneqq(-1,0);\;I_+\coloneqq\left(\sqrt{2(d-1)}-1,\sqrt{2(d-1)}\right)\;.
\end{equation}
Simply substituting shows that
\begin{equation}
    F_d(-\tfrac{1}{2})=\tfrac{1}{8}\;,
\end{equation}
independently of $d$.
Now, for every eigenvalue $\lambda\in I_0$ of $A$ there exists by the mean value theorem $\zeta\in(-1,0)$ such that
\begin{equation}\label{eq:innerevbound}
\left|\lambda+\tfrac{1}{2}\right|=\frac{\left|F_d(\lambda)-F_d(-\tfrac{1}{2})\right|}{\left|F_d'(\zeta)\right|}\le\frac{e+1/8}{2d-3}\le\frac{e}{d}\;,
\end{equation}
where we used that $F'_d(x)\le 3-2d$ on $[-1,0]$ and, in the last step, $e\ge1$, as there is no Moore graph of girth $7$ and degree $d\ge3$~\cite{HoSi60}.
For $\lambda\in I_-\cup I_+$ we only need to derive from~\eqref{eq:intervals} the trivial bound
\begin{equation}\label{eq:outerevbound}(\sqrt{2(d-1)}-1)^2<\lambda^2\;,
\end{equation}
so the eigenvalues lie in two outer clusters of width $1$ near $\pm\sqrt{2d}$ and a middle cluster of radius $e/d$ around $-\tfrac{1}{2}$.
We can also estimate how many eigenvalues lie in each cluster. Using the well-known identities for the trace and the fact that $A^2$ equals $d$ on the diagonal, we get for the number $n_{\mathrm{out}}$ of eigenvalues of $A$ in $I_-\cup I_+$ (with multiplicities) that
\begin{equation}\label{eq:boundnouter}
nd=\operatorname{tr}(A^2)=d^2+\sum_{\lambda\ne d}\lambda^2\ \ge\ d^2+n_{\mathrm{out}}\,(\sqrt{2(d-1)}-1)^2\;.
\end{equation}
This implies that only fewer than $(\tfrac{1}{2}+\tfrac{2}{\sqrt{d}})n-1$ eigenvalues may lie in the outer clusters, since otherwise
\[n_{\mathrm{out}}\cdot (\sqrt{2(d-1)}-1)^2\ge (n(\tfrac{1}{2}+\tfrac{2}{\sqrt{d}})-1)(2d-4\sqrt{d}) >nd\;.\]
So at least $(\tfrac{1}{2}-\tfrac{2}{\sqrt{d}})n$ eigenvalues are close to $-1/2$. This is summarised in Figure~\ref{fig:clusters}.\\

The key idea to arrive at a final contradiction is to shift the middle cluster of eigenvalues close to $0$ by considering $2A+I$. On the one hand this matrix is $1$ on the diagonal and $0$ or $2$ everywhere else and hence $\det(2A+I)$ is an odd integer. In particular this implies
\begin{equation}\label{eq:detodd}
|\det(2A+I)|\ge 1\;.
\end{equation}
On the other hand we can derive from~\eqref{eq:intervals} that $|2\lambda+1|<2\sqrt{2(d-1)}+1\le 3\sqrt{d}$ for $\lambda\in I_-\cup I_+$ and from~\eqref{eq:innerevbound} that $|2\lambda+1|\le\tfrac{2e}{d}$ for $\lambda\in I_0$. Thus
\begin{align*}
|\det(2A+I)|&\le (2d+1)\,(3\sqrt{d})^{n(\frac{1}{2}+\frac2{\sqrt d})}\,(2e/d)^{n(\frac12-\frac{2}{\sqrt d})}\\
&=(2d+1)\Bigl(\frac{6e}{\sqrt{d}}\Bigr)^{n/2}\Bigl(\frac{3d^{3/2}}{2e}\Bigr)^{2n/\sqrt d}\\
&\le (2d+1)\Bigl(\frac{6e}{\sqrt d}\cdot d^{8/\sqrt d}\Bigr)^{n/2}\;.
\end{align*}
Since $6e/\sqrt{d}\le 6/7$ by~\eqref{eq:exsqrtd}, and $d^{8/\sqrt{d}}\rightarrow 1$, the base is at most $0.9$ for large $d$. Therefore, as $n\ge d^3/2$, the right-hand side tends to 0, contradicting~\eqref{eq:detodd}.\\

The sketch is really wasteful in only one place. It bounds every error $|\mu_\lambda|$ by its maximum $e$. The determinant, however, only sees the product of the middle factors $|2\lambda+1|$, each of which is at most about $|\mu_\lambda|/d$. So only the geometric mean of the errors matters, and it can be shown to be of order at most $\sqrt e$. This achieves the linear bound in Theorem~\ref{thm:girth7}. The proof of the defect result, Theorem~\ref{thm:diam3}, is very similar. One defines a defect graph and obtains an equation similar to~\eqref{eq:matrixpoly}. However, unlike the excess graph, the defect graph may be a multigraph, so the averaging argument is not available, which leads to the weaker bound for the defect.
For proving Theorem~\ref{thm:girth11} in the case of girth $11$ the approach is the same, but now $F_d=F_{d,3}$ is replaced by a polynomial $F_{d,5}$ of degree $5$, and the eigenvalues lie in five clusters, near $\pm\sqrt{3d}$, $\pm\sqrt d$ and $-\frac{1}{3}$, with shares about $\frac{1}{12},\frac{1}{4},\frac{1}{3},\frac{1}{4},\frac{1}{12}$. For the final contradiction the naive choice would be to consider $3A+I$. However, for this choice the determinant is not small, as the share of the middle cluster is not large enough. Instead we consider $\det P_4(A)$ for the quartic $P_4(x)=(3x+1)(384x+371)(8x^2+4x-8d+7)+360$, which is small on the middle cluster as well as on the inner clusters close to $\pm\sqrt{d}$.

\subsection*{Organisation of the paper}
In Section~\ref{sec:prelim} we introduce the polynomials $F_{d,k}$ and their relation to non-backtracking walks, measure the deviation of a $d$-regular graph from a Moore graph by the quantities $r_k$ and $s_k$, and prove Lemma~\ref{lem:errors}, which relates the spectrum of the adjacency matrix to these quantities for every $k$. In Section~\ref{sec:reduction} we state four general theorems about arbitrary $d$-regular graphs in terms of their deviation from a Moore graph, and deduce Theorems~\ref{thm:girth7}--\ref{thm:diam5} from them. In Section~\ref{sec:seven} we prove the two statements for girth $7$ and diameter $3$ (Theorems~\ref{thm:deviation7} and~\ref{thm:deviationref7}), which sharpens the above proof sketch. In Section~\ref{sec:eleven} we prove the two statements for girth $11$ and diameter $5$ (Theorems~\ref{thm:deviation11} and~\ref{thm:deviationref11}) along the same lines, with a quartic auxiliary polynomial in place of $2x+1$. In Section~\ref{sec:remarks} we discuss the limits of the method, cages and the known constructions, and explain why the method does not apply to even girth. The few standard facts from linear algebra that we use are collected in Appendix~\ref{app:facts}.

\subsection*{Notation}
Throughout, $d\ge2$ and $k\ge1$ are integers, and
\[
m(d,k)=1+d\sum_{i=0}^{k-1}(d-1)^i
\]
is the Moore bound. All graphs are finite and simple, $n=|V(G)|$ is the number of vertices of the graph $G$ under consideration, $\dist(u,v)$ is the distance in $G$, and $A$ is the adjacency matrix of $G$. We write $\mathbf 1$ for the all-one vector, $I$ for the identity matrix and $J$ for the all-one matrix; their dimension is always $n$. The excess of a graph is denoted by $e$ and its defect by $f$.

\section{Preliminaries}\label{sec:prelim}

\subsection*{Non-backtracking walks}
A walk $v_0,v_1,\ldots,v_\ell$ is called \emph{non-backtracking} if $v_{i-1}\neq v_{i+1}$ for every $i\in\{1,\ldots,\ell-1\}$. Given a $d$-regular graph $G$, we define for every integer $\ell\ge0$ the \emph{non-backtracking walk matrix} $A^{(\ell)}$, whose $(u,v)$-entry is the number of non-backtracking walks of length $\ell$ from $u$ to $v$. Thus $A^{(0)}=I$ and $A^{(1)}=A$. Directly from the definition one obtains the well-known relations, see e.g.~\cite[Sec.~2]{ABLS07},
\begin{equation}\label{eq:nbrecursion}
A^{(2)}=A^2-dI\qquad\text{and}\qquad A^{(\ell+1)}=AA^{(\ell)}-(d-1)A^{(\ell-1)}\quad\text{for }\ell\ge2\;.
\end{equation}
Indeed, $(AA^{(\ell)})_{u,v}$ counts the walks $u,v_1,v_2,\ldots,v_{\ell+1}=v$ for which $v_1,\ldots,v_{\ell+1}$ is non-backtracking, and such a walk is non-backtracking unless $v_2=u$. For $\ell=1$ these exceptional walks are the $d$ walks $u,v_1,u$, which gives the first relation. For $\ell\ge2$ they are the walks $u,v_1,u,v_3,\ldots,v$, where $u,v_3,\ldots,v$ is a non-backtracking walk of length $\ell-1$ and $v_1$ is one of the $d-1$ neighbours of $u$ other than $v_3$, which gives the second relation.

Summing these matrices leads to the following polynomials, which have been used for Moore graphs since Hoffman and Singleton~\cite{HoSi60} and Bannai and Ito~\cite{BaIt73}; see also the survey of Miller and \v{S}ir\'a\v{n}~\cite[Sec.~3.2]{MiSi13}.

\begin{defi}\label{def:polys}
Let $F_{d,0}(x)\coloneqq1$, $F_{d,1}(x)\coloneqq x+1$, and for every $k\ge1$
\[
F_{d,k+1}(x)\coloneqq x\,F_{d,k}(x)-(d-1)F_{d,k-1}(x)\;.
\]
\end{defi}
In particular,
\begin{align}
F_{d,3}(x)&=x^3+x^2-2(d-1)x-(d-1)\;,\label{eq:F3}\\
F_{d,5}(x)&=x^5+x^4-4(d-1)x^3-3(d-1)x^2+3(d-1)^2x+(d-1)^2\;.\label{eq:F5}
\end{align}
The polynomial $F_d$ from the introduction is $F_{d,3}$.
For vertices $u,v$ of a $d$-regular graph we denote by $\nb_k(u,v)$ the number of non-backtracking walks of length at most $k$ from $u$ to $v$, where the walk of length $0$ counts if $u=v$.

\begin{lem}\label{lem:polynomials}
Let $G$ be a $d$-regular graph with adjacency matrix $A$ and let $k\ge0$. Then
\[
F_{d,k}(A)=\sum_{i=0}^k A^{(i)}=\bigl(\nb_k(u,v)\bigr)_{u,v}\;.
\]
In particular, all row sums of $F_{d,k}(A)$ are equal to $F_{d,k}(d)=m(d,k)$.
\end{lem}
\begin{proof}
The second equality is the definition of $\nb_k$. For the first, let $S_k\coloneqq\sum_{i=0}^kA^{(i)}$. Then $S_0=I=F_{d,0}(A)$ and $S_1=I+A=F_{d,1}(A)$, so it suffices to show that $S_{k+1}-AS_k+(d-1)S_{k-1}=0$ for every $k\ge1$. Collecting the terms with the same upper index, this expression equals
\[
A^{(0)}+\bigl(A^{(1)}-AA^{(0)}\bigr)+\bigl(A^{(2)}-AA^{(1)}+(d-1)A^{(0)}\bigr)+\sum_{i=3}^{k+1}\bigl(A^{(i)}-AA^{(i-1)}+(d-1)A^{(i-2)}\bigr)\;,
\]
which by~\eqref{eq:nbrecursion} is $I+0-I+0=0$. Finally, a vertex is the starting point of exactly $d(d-1)^{i-1}$ non-backtracking walks of length $i\ge1$, so every row sum of $S_k$ is $1+d\sum_{i=1}^{k}(d-1)^{i-1}=m(d,k)$. Since $A\mathbf1=d\mathbf1$, the row sums of $F_{d,k}(A)$ are also equal to $F_{d,k}(d)$.
\end{proof}

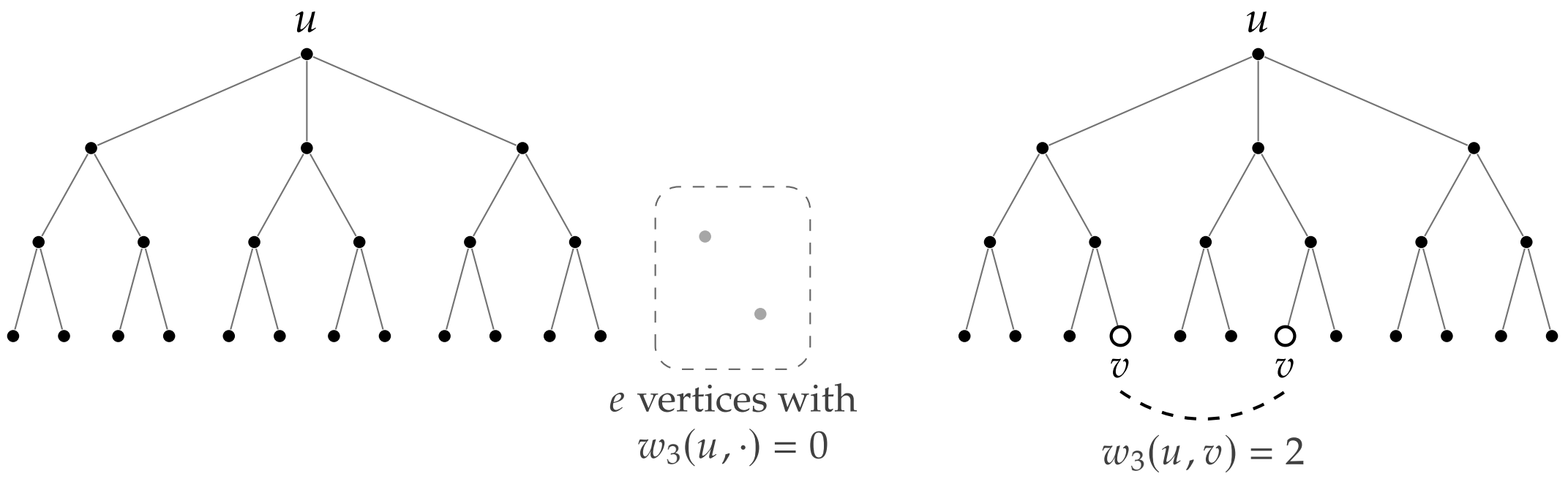
\begin{figure}[ht]
\centering
\begin{tikzpicture}[
  v/.style={circle, fill=black, inner sep=1.1pt},
  w/.style={circle, fill=black!35, inner sep=1.1pt},
  level 1/.style={sibling distance=1.95cm, level distance=0.85cm},
  level 2/.style={sibling distance=0.95cm, level distance=0.85cm},
  level 3/.style={sibling distance=0.46cm, level distance=0.85cm},
  edge from parent/.style={draw=black!55}]
\begin{scope}
\node[v, label=above:{$u$}] {}
  child {node[v]{} child{node[v]{} child{node[v]{}} child{node[v]{}}} child{node[v]{} child{node[v]{}} child{node[v]{}}}}
  child {node[v]{} child{node[v]{} child{node[v]{}} child{node[v]{}}} child{node[v]{} child{node[v]{}} child{node[v]{}}}}
  child {node[v]{} child{node[v]{} child{node[v]{}} child{node[v]{}}} child{node[v]{} child{node[v]{}} child{node[v]{}}}};
\draw[rounded corners=6pt, dashed, black!60] (3.15,-2.85) rectangle (4.55,-1.2);
\node[w] at (3.6,-1.65) {};
\node[w] at (4.1,-2.35) {};
\node[align=center, font=\small, black!75] at (3.85,-3.35) {$e$ vertices with\\ $\nb_3(u,\cdot)=0$};
\node[font=\small] at (0.6,-4.1) {(a) girth $\ge7$, excess $e$};
\end{scope}
\begin{scope}[xshift=8.6cm]
\node[v, label=above:{$u$}] {}
  child {node[v]{} child{node[v]{} child{node[v]{}} child{node[v]{}}} child{node[v]{} child{node[v]{}} child{node[v] (y1) {}}}}
  child {node[v]{} child{node[v]{} child{node[v]{}} child{node[v]{}}} child{node[v]{} child{node[v] (y2) {}} child{node[v]{}}}}
  child {node[v]{} child{node[v]{} child{node[v]{}} child{node[v]{}}} child{node[v]{} child{node[v]{}} child{node[v]{}}}};
\draw[thick, fill=white] (y1) circle (2.4pt) (y2) circle (2.4pt);
\node[below=2pt, font=\small] at (y1) {$v$};
\node[below=2pt, font=\small] at (y2) {$v$};
\draw[dashed, thick] ($(y1)+(0,-0.5)$) to[bend right=35] node[below, font=\small, black!75]{$\nb_3(u,v)=2$} ($(y2)+(0,-0.5)$);
\node[font=\small] at (0,-4.1) {(b) diameter $\le3$, defect $f$};
\end{scope}
\end{tikzpicture}
\caption{The non-backtracking tree of depth $3$ rooted at $u$, drawn for $d=3$; $\nb_3(u,v)$ is the number of its nodes that are copies of $v$. (a) In a graph of girth at least $7$ the nodes are distinct vertices, and the $e$ vertices at distance at least $4$ from $u$ do not appear. (b) In a graph of diameter at most $3$ every vertex appears, and in total $f$ nodes are repeated copies.}
\label{fig:trees}
\end{figure}

Equivalently, $\nb_k(u,v)$ is the number of copies of $v$ in the non-backtracking tree of depth $k$ rooted at $u$, whose nodes are the non-backtracking walks of length at most $k$ starting at $u$, each labelled by its end vertex. This tree has $m(d,k)$ nodes, so the Moore bound can be read off from the tree of a single vertex, and the two extremal problems correspond to the two ways in which the tree can fail to contain every vertex exactly once; see Figure~\ref{fig:trees}.

If $G$ is a $d$-regular graph of girth at least $2k+1$ and excess $e$, then the nodes of the tree are distinct vertices. Hence $(F_{d,k}(A))_{u,v}$ equals $1$ if $\dist(u,v)\le k$ and $0$ otherwise, and the zeros in the row of $u$ are exactly the $n-m(d,k)=e$ vertices that the tree misses. For $k=3$ this is the identity~\eqref{eq:matrixpoly}.
If $G$ is a $d$-regular graph of diameter at most $k$ and defect $f$, then every vertex appears in the tree, so $\nb_k(u,v)\ge1$ for all $u,v$, and in total $m(d,k)-n=f$ nodes of the tree are repeated copies.

\subsection*{Deviation from a Moore graph}
For a $d$-regular graph $G$ we define the symmetric integer matrix
\begin{equation}\label{eq:defCk}
C_k\coloneqq F_{d,k}(A)-J=\bigl(\nb_k(u,v)-1\bigr)_{u,v}\;.
\end{equation}
We measure how far $G$ is from a Moore graph by the largest absolute row sum of $C_k$ and by its normalised squared Frobenius norm,
\[
r_k(G)\coloneqq\max_{u\in V(G)}\ \sum_{v\in V(G)}\bigl|\nb_k(u,v)-1\bigr|
\qquad\text{and}\qquad
s_k(G)\coloneqq\frac1n\sum_{u,v\in V(G)}\bigl(\nb_k(u,v)-1\bigr)^2=\frac1n\operatorname{tr}\bigl(C_k^2\bigr)\;.
\]
So $r_k(G)$ is the largest number of vertices that the non-backtracking tree of depth $k$ of a vertex misses or repeats, counted with multiplicity. We call it the \emph{deviation} of $G$. A $d$-regular graph satisfies $r_k(G)=0$ if and only if it is a Moore graph of diameter $k$, that is, of girth $2k+1$. Since $\sum_v(\nb_k(u,v)-1)^2\le\max_v|\nb_k(u,v)-1|\cdot\sum_v|\nb_k(u,v)-1|$ for every $u$, we have
\begin{equation}\label{eq:sbound}
s_k(G)\ \le\ c_k(G)\,r_k(G)\ \le\ r_k(G)^2,\qquad\text{where }c_k(G)\coloneqq\max_{u,v}\bigl|\nb_k(u,v)-1\bigr|\;.
\end{equation}

For a graph $G$ of girth at least $2k+1$ we define its \emph{excess graph} $E$ as the graph on $V(G)$ in which two vertices are adjacent if they have distance at least $k+1$ in $G$; for $k=3$ this is the excess graph from the introduction.
We only consider graphs whose excess or defect is less than $(m(d,k)-1)/d=\sum_{i=0}^{k-1}(d-1)^i$. For such graphs regularity does not have to be imposed, since it follows from the Moore bound, and the matrix $C_k$ has a simple structure.

\begin{prop}\label{pro:regular}
Let $d\ge2$ and $k\ge1$.
\begin{enumerate}[label=\textup{(\alph*)}]
\item Let $G$ be a graph of minimum degree $d$, girth at least $2k+1$ and excess $0\le e<\sum_{i=0}^{k-1}(d-1)^i$. Then $G$ is $d$-regular, and $\nb_k(u,v)=1$ if $\dist(u,v)\le k$ and $\nb_k(u,v)=0$ otherwise. Equivalently, $F_{d,k}(A)=J-B$, where $B$ is the adjacency matrix of the excess graph $E$. Moreover, $E$ is $e$-regular and $r_k(G)=s_k(G)=e$.
\item Let $G$ be a graph of maximum degree $d$, diameter at most $k$ and defect $0\le f<\sum_{i=0}^{k-1}(d-1)^i$. Then $G$ is $d$-regular, and $C_k$ has non-negative integer entries and all row sums equal to $f$. In particular $r_k(G)=f$. If moreover $\nb_k(u,v)\le2$ for all $u,v$, then $s_k(G)=f$.
\end{enumerate}
\end{prop}
\begin{proof}
(a) Since $G$ has no cycle of length at most $2k$, every non-backtracking walk of length at most $k$ is a path, and two different such walks starting at the same vertex end at different vertices. Hence, for every vertex $v$, the breadth-first search tree of depth $k$ rooted at $v$ consists of distinct vertices. In it, $v$ has $\deg(v)$ children and every other vertex at distance less than $k$ from $v$ has at least $d-1$ children. So at least $1+\deg(v)\sum_{i=0}^{k-1}(d-1)^i$ vertices are at distance at most $k$ from $v$. If $\deg(v)\ge d+1$, this is at least $m(d,k)+\sum_{i=0}^{k-1}(d-1)^i>m(d,k)+e=n$, which is impossible. So $G$ is $d$-regular. By the first observation, $\nb_k(u,v)=1$ if $\dist(u,v)\le k$ and $\nb_k(u,v)=0$ otherwise, which together with Lemma~\ref{lem:polynomials} is the identity $F_{d,k}(A)=J-B$. By Lemma~\ref{lem:polynomials}, every row of $F_{d,k}(A)$ contains exactly $m(d,k)$ many ones, hence exactly $n-m(d,k)=e$ zeros. So $E$ is $e$-regular, $C_k=-B$, and $r_k(G)=s_k(G)=e$.

(b) Suppose that some vertex $v$ has $\deg(v)\le d-1$. Since all degrees are at most $d$, at most
\[
1+(d-1)\sum_{i=0}^{k-1}(d-1)^i=m(d,k)-\sum_{i=0}^{k-1}(d-1)^i<m(d,k)-f=n
\]
vertices are at distance at most $k$ from $v$, contradicting diameter at most $k$. So $G$ is $d$-regular. A shortest path is a non-backtracking walk, so $\nb_k(u,v)\ge1$ for all $u,v$, that is, $C_k$ has non-negative integer entries. By Lemma~\ref{lem:polynomials} its row sums are $m(d,k)-n=f$, hence $r_k(G)=f$. If moreover $\nb_k(u,v)\le2$ for all $u,v$, then all entries of $C_k$ lie in $\{0,1\}$ and are equal to their squares, so $s_k(G)=\frac1n\cdot nf=f$.
\end{proof}

The following lemma is the starting point for our proof in both cases. It turns Lemma~\ref{lem:polynomials} into a statement about the spectrum of $A$ and relates the eigenvalues to the deviation from the Moore bound.

\begin{lem}\label{lem:errors}
Let $G$ be a $d$-regular graph on $n$ vertices with adjacency matrix $A$, let $k\ge1$ and let $r\ge r_k(G)$.
\begin{enumerate}[label=\textup{(\alph*)}]
\item $|n-m(d,k)|\le r$.\label{itm:rboundserror}
\item $A\mathbf1=d\mathbf1$, and $\mathbf1^\perp$ (the subspace orthogonal to $\mathbf 1$) has an orthonormal basis of eigenvectors of $A$; let $\sigma^*$ be the multiset of their $n-1$ eigenvalues. If $b$ is one of them, with eigenvalue $\lambda$, then $C_kb=\mu_\lambda b$ for $\mu_\lambda\coloneqq F_{d,k}(\lambda)$; in particular $|\mu_\lambda|\le r$.\label{itm:eigenvofC}
\item\label{itm:eigenvands(G)} We have
\begin{equation}\label{eq:moments}
\sum_{\lambda\in\sigma^*}\lambda=-d,\qquad\sum_{\lambda\in\sigma^*}\lambda^2=nd-d^2\qquad\text{and}\qquad\sum_{\lambda\in\sigma^*}\mu_\lambda^2=n\,s_k(G)-\bigl(n-m(d,k)\bigr)^2 .
\end{equation}
\end{enumerate}
\end{lem}
\begin{proof}
To show~\ref{itm:rboundserror} note that by Lemma~\ref{lem:polynomials}, every row of $C_k$ sums to $m(d,k)-n$. Hence $|m(d,k)-n|\le\sum_v|(C_k)_{u,v}|\le r_k(G)\le r$.
For~\ref{itm:eigenvofC} we have $A\mathbf1=d\mathbf1$, since $G$ is $d$-regular. Since $A$ is a real symmetric matrix, $\mathbf1/\sqrt n$ can be extended to an orthonormal basis of eigenvectors (Fact~\ref{fact:spectral}), and the other basis vectors form a basis of $\mathbf1^\perp$. If $Ab=\lambda b$ with $b\perp\mathbf1$, then $Jb=0$, so by~\eqref{eq:defCk} and Fact~\ref{fact:trace}
\[
C_kb=\bigl(F_{d,k}(A)-J\bigr)b=F_{d,k}(A)b=F_{d,k}(\lambda)b\;.
\]
The matrix $C_k$ is symmetric with largest absolute row sum $r_k(G)\le r$, and hence by Fact~\ref{fact:rowsum} all its eigenvalues lie in $[-r,r]$.
Lastly, for~\ref{itm:eigenvands(G)} by Fact~\ref{fact:trace}, $\sum_{\lambda\in\sigma^*}\lambda=\operatorname{tr}A-d=-d$ and $\sum_{\lambda\in\sigma^*}\lambda^2=\operatorname{tr}A^2-d^2=nd-d^2$. Moreover, $C_k\mathbf1=(m(d,k)-n)\mathbf1$, so by~\ref{itm:eigenvofC} the eigenvalues of $C_k^2$ are $(m(d,k)-n)^2$ and the numbers $\mu_\lambda^2$, while $\operatorname{tr}C_k^2=n\,s_k(G)$.
\end{proof}

We call the elements of $\sigma^*$ the \emph{non-trivial eigenvalues} of $A$, and $\mu_\lambda$ the \emph{error} of $\lambda$. In the situation of Proposition~\ref{pro:regular}(a) we have $C_k=-B$, so $-\mu_\lambda$ is an eigenvalue of the excess graph; for $k=3$, the error $\mu_\lambda=F_d(\lambda)$ is the error from the proof sketch. In any case, every non-trivial eigenvalue is a root of $F_{d,k}(x)=t$ for some $t\in[-r,r]$. For $k\in\{3,5\}$ and $r$ small compared with $d$, this confines the non-trivial eigenvalues to a few short intervals, and the rest of the paper exploits this.

\section{Reduction to the deviation from a Moore graph}\label{sec:reduction}

We deduce Theorems~\ref{thm:girth7}--\ref{thm:diam5} from four statements about arbitrary $d$-regular graphs, in which the excess and the defect are replaced by the deviation $r_k(G)$ from a Moore graph, or by the mean square deviation $s_k(G)$, defined in Section~\ref{sec:prelim}. For girth $7$ and diameter $3$ we have the following statements, which we prove in Section~\ref{sec:seven}.

\begin{thm}\label{thm:deviation7}
Every $d$-regular graph $G$ satisfies $d<\bigl(\sqrt{8}\,r_3(G)+2\bigr)^2$.
\end{thm}

To achieve the linear bound for the girth case we establish a stronger conditional relation.

\begin{thm}\label{thm:deviationref7}
Let $G$ be a $d$-regular graph such that $C_3$ or $-C_3$ has all entries in $\{0,1\}$, that is, either $\nb_3(u,v)\le1$ for all $u,v$, or $1\le\nb_3(u,v)\le2$ for all $u,v$. Then $d<\bigl(\sqrt{8\,r_3(G)}+2\bigr)^2$.
\end{thm}

For girth $11$ and diameter $5$ we show the following two statements, with absolute constants that we do not track. We prove them in Section~\ref{sec:eleven}.

\begin{thm}\label{thm:deviation11}
There is an absolute constant $C$ such that every $d$-regular graph $G$ satisfies $d<C\,(r_5(G)+1)^{10/3}$.
\end{thm}

\begin{thm}\label{thm:deviationref11}
There is an absolute constant $C$ such that every $d$-regular graph $G$ with $d\ge C\,(r_5(G)+1)$ satisfies $d<C\,(s_5(G)+1)^{5/3}$.
\end{thm}

We now deduce our main theorems. In all four proofs, Proposition~\ref{pro:regular} turns the excess or the defect into the deviation.

\begin{proof}[Proof of Theorem~\ref{thm:girth7}]
Let $G$ have minimum degree $d$, girth at least $7$ and excess $e$. If $e\ge 1+(d-1)+(d-1)^2$, then $(\sqrt{8e}+2)^2>8e>d$. So we may assume otherwise and by Proposition~\ref{pro:regular}(a), $G$ is $d$-regular and $-C_3=B$ is the adjacency matrix of the excess graph $E$ and has all entries in $\{0,1\}$ and $r_3(G)=e$.
Hence Theorem~\ref{thm:deviationref7} gives $d<(\sqrt{8e}+2)^2$. This proves the second statement, for every $d\ge3$. For $d\ge4$ we have $\sqrt{d}-2\ge0$, so $\sqrt d<\sqrt{8e}+2$ gives $e>(\sqrt d-2)^2/8$, which is the first statement. For $d=3$ there is no Moore graph of degree $3$ and girth $7$~\cite{HoSi60}, so $e\ge1>(\sqrt{3}-2)^2/8$.\end{proof}

\begin{proof}[Proof of Theorem~\ref{thm:diam3}]
Let $G$ have maximum degree $d$, diameter at most $3$ and defect $f$. If $f\ge1+(d-1)+(d-1)^2$, then $(\sqrt8\,f+2)^2>8f^2>d$. Otherwise, by Proposition~\ref{pro:regular}(b), $G$ is $d$-regular with $r_3(G)=f$, and Theorem~\ref{thm:deviation7} gives $d<(\sqrt8\,f+2)^2$. For $d\ge4$, solving for $f$ gives $f>(\sqrt d-2)/\sqrt8$, which is the first statement. For $d=3$ the first statement says that $G$ has fewer than $m(d,3)+(2-\sqrt3)/\sqrt8$ vertices, which holds by the Moore bound.
\end{proof}

\begin{proof}[Proof of Theorems~\ref{thm:girth11} and~\ref{thm:diam5}]
Let $C\ge1$ be a constant for which Theorems~\ref{thm:deviation11} and~\ref{thm:deviationref11} hold, and let $d\ge3$.

Let $G$ have minimum degree $d$, girth at least $11$ and excess $e$. If $e\ge1+(d-1)+\dots+(d-1)^4$, then $e\ge d$ and $d<C(e+1)^{5/3}$ trivially. Otherwise $G$ is $d$-regular with $r_5(G)=s_5(G)=e$ by Proposition~\ref{pro:regular}(a). If $d<C(e+1)$, then $d<C(e+1)^{5/3}$, and otherwise Theorem~\ref{thm:deviationref11} gives $d<C(e+1)^{5/3}$. So in all cases $e+1>(d/C)^{3/5}$. Moreover $e\ge1$, as there is no Moore graph of degree $d\ge3$ and girth $11$~\cite{BaIt73,Dam73}. Hence $2e\ge e+1>(d/C)^{3/5}$, which is Theorem~\ref{thm:girth11} with $c=\frac{1}{2}C^{-3/5}$.

Let $G$ have maximum degree $d$, diameter at most $5$ and defect $f$. If $f\ge1+(d-1)+\dots+(d-1)^4$, then $d<C(f+1)^{10/3}$ trivially. Otherwise $G$ is $d$-regular with $r_5(G)=f$ by Proposition~\ref{pro:regular}(b), and Theorem~\ref{thm:deviation11} gives $d<C(f+1)^{10/3}$. As before $f\ge1$~\cite{BaIt73,Dam73}, so $2f\ge f+1>(d/C)^{3/10}$, which is Theorem~\ref{thm:diam5} with $c=\frac{1}{2}C^{-3/10}$.
\end{proof}

\section{Girth seven and diameter three}\label{sec:seven}

In this section we prove Theorems~\ref{thm:deviation7} and~\ref{thm:deviationref7}.

\subsection*{Setup}
Throughout this section, $G$ is a $d$-regular graph on $n$ vertices with adjacency matrix $A$, and $k=3$. We write $z\coloneqq d-1$ and abbreviate
\[ F\coloneqq F_{d,3}=x^3+x^2-2zx-z,\qquad C\coloneqq C_3=F(A)-J,\qquad m\coloneqq m(d,3)=d^3-d^2+d+1, \]
as well as $r(G)\coloneqq r_3(G)$ and $s(G)\coloneqq s_3(G)$. Moreover, $r$ is a real number such that
\begin{equation}\label{eq:zlarge}
r\ge r(G)\qquad\text{and}\qquad z\ge\max(2r+2,3)\;.
\end{equation}
At the beginning of the proofs of Theorems~\ref{thm:deviation7} and~\ref{thm:deviationref7} we check that these assumptions hold for $r=r(G)$. As in Section~\ref{sec:prelim}, $\sigma^*$ is the multiset of the $n-1$ non-trivial eigenvalues of $A$, and $\mu_\lambda=F(\lambda)$ is the error of $\lambda\in\sigma^*$. By Lemma~\ref{lem:errors} with $k=3$,
\begin{gather}
|n-m|\le r\qquad\text{and}\qquad|\mu_\lambda|\le r\quad\text{for all }\lambda\in\sigma^*,\label{eq:input1}\\
\sum_{\lambda\in\sigma^*}\lambda=-d,\qquad\sum_{\lambda\in\sigma^*}\lambda^2=nd-d^2\qquad\text{and}\qquad\sum_{\lambda\in\sigma^*}\mu_\lambda^2\le n\,s(G)\;.\label{eq:input2}
\end{gather}

\subsection*{Outline}
The proof refines the proof of~\eqref{eq:weakerbound} in the introduction in five steps.
\begin{enumerate}[label=\textup{(\arabic*)}]
\item By~\eqref{eq:input1}, every non-trivial eigenvalue $\lambda$ of $A$ is a root of $F(x)=\mu_\lambda$ with $|\mu_\lambda|\le r$.
\item Hence $\lambda$ lies in one of three short clusters: a middle one near $-\frac{1}{2}$, in which $|2\lambda+1|\le(1+8|\mu_\lambda|)/(8z-4)$, and two outer ones near $\pm\sqrt{2z}$ (Lemma~\ref{lem:clusters}). We write $\sigma_{\mathrm{mid}}$ and $\sigma_{\mathrm{out}}$ for the multisets of middle and outer non-trivial eigenvalues, and $n_{\mathrm{mid}}$ and $n_{\mathrm{out}}$ for their sizes.
\item Analysing the second moment $\operatorname{tr}A^2=nd$ forces at least about half of the eigenvalues into the middle cluster (Lemma~\ref{lem:counting}).
\item Since $\det(2A+I)$ is an odd integer, $\prod_{\lambda\in\sigma^*}|2\lambda+1|\ge1/(2d+1)$. As each outer factor is only about $\sqrt{8d}$, the middle factors satisfy $\prod_{\lambda\in\sigma_{\mathrm{mid}}}|2\lambda+1|\gtrsim(8d)^{-n_{\mathrm{out}}/2}$ (Lemma~\ref{lem:det}).
\item By (2), each middle factor is at most about $(1+8|\mu_\lambda|)/(8d)$. As $n_{\mathrm{out}}$ is at most about $n_{\mathrm{mid}}$ by (3), step (4) forces $1+8|\mu_\lambda|\gtrsim\sqrt{8d}$ on geometric average over the middle cluster. Hence the root mean square $\rho$ of the middle errors satisfies $1+8\rho\ge\Phi$ for an explicit $\Phi=\sqrt{8d}-O(\log d/\sqrt d)$ (Proposition~\ref{prop:key} and Lemma~\ref{lem:phi}).
\end{enumerate}
Since $\rho\le r(G)$, step (5) gives Theorem~\ref{thm:deviation7}. For Theorem~\ref{thm:deviationref7} we show $\rho^2\lesssim r(G)$. If $C$ or $-C$ is a $0/1$ matrix, the mean square of all errors is only about $r(G)$, and at most about half of it sits on the middle cluster (Lemma~\ref{lem:halfmass}).

\subsection*{Three clusters}
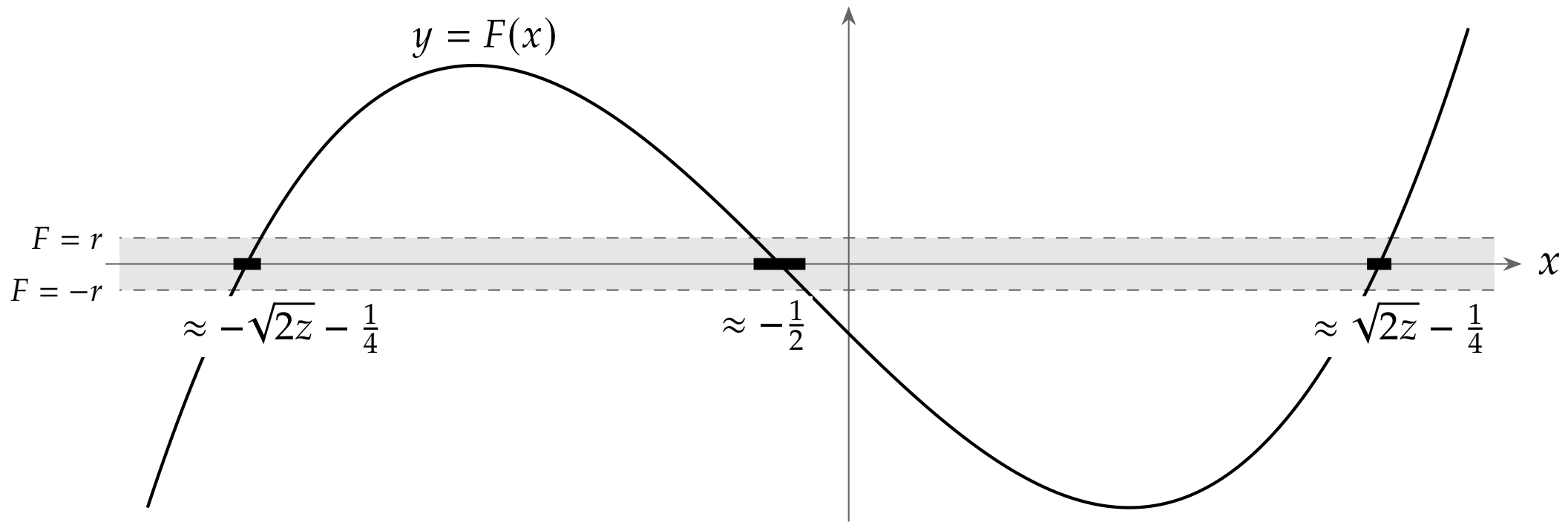
\begin{figure}[ht]
\centering
\begin{tikzpicture}[x=1.25cm, y=0.62cm]
\fill[black!10] (-5.2,-0.375) rectangle (4.6,0.375);
\draw[black!60, thin, dashed] (4.6,0.375) -- (-5.2,0.375) node[left, font=\scriptsize, black] {$F=r$};
\draw[black!60, thin, dashed] (4.6,-0.375) -- (-5.2,-0.375) node[left, font=\scriptsize, black] {$F=-r$};
\draw[-{Stealth}, black!60] (-5.3,0) -- (4.8,0) node[right, black] {$x$};
\draw[-{Stealth}, black!60] (0,-3.7) -- (0,3.7);
\draw[thick, smooth, samples=160, domain=-5.0:4.42] plot (\x, {((\x)*(\x)*(\x)+(\x)*(\x)-16*(\x)-8)/8});
\draw[line width=3pt, black] (-4.3872,0) -- (-4.1913,0);
\draw[line width=3pt, black] (-0.6782,0) -- (-0.3084,0);
\draw[line width=3pt, black] (3.6956,0) -- (3.8695,0);
\node[below=8pt, font=\small, fill=white, inner sep=1pt] at (-4.04,0) {$\approx-\sqrt{2z}-\frac{1}{4}$};
\node[below=8pt, font=\small, fill=white, inner sep=1pt] at (-0.6,0) {$\approx-\frac{1}{2}$};
\node[below=8pt, font=\small, fill=white, inner sep=1pt] at (3.93,0) {$\approx\sqrt{2z}-\frac{1}{4}$};
\node[font=\small] at (-2.6,3.25) {$y=F(x)$};
\end{tikzpicture}
\caption{The polynomial $F=F_{9,3}$ together with the strip $|F(x)|\le r$ for $r=3$. (The vertical axis is compressed.) The preimage $F^{-1}([-r,r])$ consists of three short intervals (thick), which contain all non-trivial eigenvalues.}
\label{fig:F}
\end{figure}

The next lemma locates the real roots of $F(x)=t$ for $|t|\le r$, and hence the non-trivial eigenvalues, as shown in Figure~\ref{fig:F}. Its two types correspond to the interval $I_0$ and to $I_-\cup I_+$ in~\eqref{eq:intervals}. For the rest of this section we set
\[
Q(x)\coloneqq x^2+\frac{x}{2},\qquad\eps\coloneqq\frac{8r+1}{8z-4}\qquad\text{and}\qquad\tau\coloneqq\frac{r}{\sqrt{2z}}\;.
\]

\begin{lem}[Cluster Lemma]\label{lem:clusters}
Let $|t|\le r$. Every real root $\theta$ of $F(x)=t$ satisfies exactly one of the following:
\begin{enumerate}
\labelleditem{(M)}{itm:middlecl} $|2\theta+1|\le\dfrac{1+8|t|}{8z-4}\le\eps$ and $|Q(\theta)|\le\dfrac{\eps}{2}$ \quad (middle cluster);
\labelleditem{(O)}{itm:outercl} $\sqrt{2z}-1\le|\theta|\le\sqrt{2z}+1$ and $|Q(\theta)-2z|\le\dfrac{1}{2}+\dfrac{|t|}{\sqrt{2z}}\le\dfrac{1}{2}+\tau$ \quad (outer clusters).
\end{enumerate}
\end{lem}
\begin{proof}
By~\eqref{eq:zlarge} we have $|t|\le r<z$ and $\eps<1$. In particular, the two types are mutually exclusive: type~\ref{itm:middlecl} gives $|2\theta+1|\le\eps<1$ and hence $|\theta|<1$, whereas type~\ref{itm:outercl} gives $|\theta|\ge\sqrt{2z}-1>1$ as $z\ge3$. So it suffices to show that $\theta$ has one of the two types. We use the two easily checked identities
\begin{equation}\label{eq:Fids}
F(x)=(x^2-2z)(x+1)+z\qquad\text{and}\qquad F(x)-\frac{1}{8}=\Bigl(x+\frac{1}{2}\Bigr)\Bigl(Q(x)-2z-\frac{1}{4}\Bigr)\;.
\end{equation}
By the first identity,
\[
\begin{array}{c|cccccc}
x & -\sqrt{2z}-1 & -\sqrt{2z} & -1 & 0 & \sqrt{2z}-1 & \sqrt{2z}\\[2pt]
\hline
F(x)\rule{0pt}{12pt} & -3z-\sqrt{2z} & z & z & -z & \sqrt{2z}-3z & z
\end{array}\;.
\]
As $|t|<z$ and $\sqrt{2z}\le2z$, the cubic polynomial $F(x)-t$ is negative, positive, positive, negative, negative and positive at these six points. So it has a root in each of the intervals $I_-$, $I_0$ and $I_+$ from~\eqref{eq:intervals}, and these are all its real roots. Moreover, since $F(\theta)=t$, the second identity gives
\begin{equation}\label{eq:thetaid}
\Bigl(\theta+\frac{1}{2}\Bigr)\Bigl(Q(\theta)-2z-\frac{1}{4}\Bigr)=t-\frac{1}{8}\;.
\end{equation}

If $\theta\in I_0=(-1,0)$, then $Q(\theta)=\theta(\theta+\frac{1}{2})<\frac{1}{2}$, so the second factor in~\eqref{eq:thetaid} has absolute value $2z+\frac{1}{4}-Q(\theta)>2z-1$. Hence
\[
|2\theta+1|\ \eqByRef{eq:thetaid}\ \frac{2|t-\frac{1}{8}|}{2z+\frac{1}{4}-Q(\theta)}\ \le\ \frac{2|t|+\frac{1}{4}}{2z-1}\ =\ \frac{1+8|t|}{8z-4}\ \le\ \eps\;.
\]
With $y\coloneqq\theta+\frac{1}{2}$ we have $Q(\theta)=y^2-\frac{y}{2}$ and $|y|\le\frac{\eps}{2}$, so $|Q(\theta)|\le y^2+\frac{|y|}{2}\le\frac{\eps^2}{4}+\frac{\eps}{4}\le\frac{\eps}{2}$, as $\eps<1$. So $\theta$ is of type~\ref{itm:middlecl}.

If $\theta\in I_-\cup I_+$, then $\sqrt{2z}-1<|\theta|<\sqrt{2z}+1$ and $|\theta+\frac{1}{2}|>\sqrt{2z}-\frac{1}{2}$. Hence
\[
\Bigl|Q(\theta)-2z-\frac{1}{4}\Bigr|\ \eqByRef{eq:thetaid}\ \frac{|t-\frac{1}{8}|}{|\theta+\frac{1}{2}|}\ \le\ \frac{|t|+\frac{1}{8}}{\sqrt{2z}-\frac{1}{2}}\ \le\ \frac{1}{4}+\frac{|t|}{\sqrt{2z}}\;,
\]
where the last inequality can be shown by multiplying both sides by $\sqrt{2z}\,(\sqrt{2z}-\frac{1}{2})$ leading to $|t|+\sqrt{z/2}\le z$, which holds as $|t|\le\frac{z}{2}-1$ and $\sqrt{z/2}\le\frac{z}{2}$ by~\eqref{eq:zlarge}. Thus $|Q(\theta)-2z|\le\frac{1}{2}+\frac{|t|}{\sqrt{2z}}$, and $\theta$ is of type~\ref{itm:outercl}.
\end{proof}

\begin{rem}
The bound on $Q(\theta)$ in~\ref{itm:outercl} is a refined version of the bound $(\sqrt{2z}-1)^2<\theta^2$ in~\eqref{eq:outerevbound}. Its advantage is symmetry: by~\eqref{eq:Fids}, $Q(\theta)\approx2z+\frac{1}{4}$ on both outer clusters, whereas $\theta^2\approx2z-\sqrt{z/2}$ on the right one and $\theta^2\approx2z+\sqrt{z/2}$ on the left one.
\end{rem}

\subsection*{Half of the spectrum is close to \texorpdfstring{$-\frac{1}{2}$}{-1/2}}
By~\eqref{eq:input1}, Lemma~\ref{lem:clusters} applies to every non-trivial eigenvalue $\lambda$, with $\theta=\lambda$ and $t=\mu_\lambda$. We call $\lambda$ \emph{middle} or \emph{outer} according to its type. We denote by $\sigma_{\mathrm{mid}}$ and $\sigma_{\mathrm{out}}$ the sub-multisets of middle and of outer eigenvalues in $\sigma^*$, and by $n_{\mathrm{mid}}$ and $n_{\mathrm{out}}$ their sizes, so that $n_{\mathrm{mid}}+n_{\mathrm{out}}=n-1$.

\begin{lem}[Cluster Size Lemma]\label{lem:counting}
Let $s\ge s(G)$ and
\[
\eta\coloneqq\frac{\frac{5}{2}+\eps+2\sqrt{s/(2z)}}{2d-\frac{5}{2}}\;.
\]
Then $n_{\mathrm{out}}\le\frac{n}{2}(1+\eta)$. Suppose moreover that $\eps\le\frac{1}{8}$ and $s\le\frac{d}{8}$. If $d\ge24$, then $\eta<\frac{1}{10}$, and if even $d\ge100$, then $\eta\le1.6/d$.
\end{lem}
\begin{proof}
By the Cauchy--Schwarz inequality,
\[
\sum_{\lambda\in\sigma_{\mathrm{out}}}|\mu_\lambda|\ \le\ \Bigl(n\sum_{\lambda\in\sigma^*}\mu_\lambda^2\Bigr)^{1/2}\ \leByRef{eq:input2}\ \bigl(n\cdot ns\bigr)^{1/2}\ =\ n\sqrt{s}\;,
\]
using $s(G)\le s$ in the second step.
By Lemma~\ref{lem:clusters}, $Q(\lambda)\ge-\frac{\eps}{2}$ for middle and $Q(\lambda)\ge2z-\frac{1}{2}-\frac{|\mu_\lambda|}{\sqrt{2z}}$ for outer $\lambda$. Hence
\[
\Bigl(2z-\frac{1}{2}\Bigr)n_{\mathrm{out}}-\frac{n\sqrt{s}}{\sqrt{2z}}-\frac{\eps}{2}\,n\ \le\ \sum_{\lambda\in\sigma^*} Q(\lambda)\
=\ \sum_{\lambda\in\sigma^*}\Bigl(\lambda^2+\frac{\lambda}{2}\Bigr)\
\eqByRef{eq:input2}\ nd-d^2-\frac{d}{2}\ <\ nd\;.
\]
Solving for $n_{\mathrm{out}}$ and using $2z-\frac{1}{2}=2d-\frac{5}{2}$ gives
\[
n_{\mathrm{out}}\ <\ \frac{n}{2}\cdot\frac{2d+\eps+2\sqrt{s/(2z)}}{2d-\frac{5}{2}}\ =\ \frac{n}{2}(1+\eta)\;.
\]
Under the additional assumptions, $2\sqrt{s/(2z)}\le\frac{1}{2}\sqrt{d/z}$, which is less than $0.511$ for $d\ge24$. So $\eta\le3.14/(2d-\frac{5}{2})$, which is less than $\frac{1}{10}$ for $d\ge24$ and at most $1.6/d$ for $d\ge100$.
\end{proof}

Since $\eta=O(1/d)$ for $s=O(d)$, and $n\approx d^3$, Lemma~\ref{lem:counting} derives from the second moment that at most about $\frac{n}{2}$ non-trivial eigenvalues are outer, so that at least about $\frac{n}{2}$ of them lie within $\frac{\eps}{2}\approx\frac{8r+1}{16d}$ of $-\frac{1}{2}$. This is the precise form of~\eqref{eq:boundnouter}.

\subsection*{Bounding the determinant}
We bound $|\det(2A+I)|$ from below by $1$, using parity, and from above in terms of the middle eigenvalues. This gives a lower bound on the product of the numbers $|2\lambda+1|$ over the middle eigenvalues, that is, on how close these eigenvalues can be to $-\frac{1}{2}$.

\begin{lem}[Determinant Lemma]\label{lem:det}
We have
\begin{equation}\label{eq:det}
\prod_{\lambda\in\sigma_{\mathrm{mid}}}\frac{1}{|2\lambda+1|}\ \le\ (2d+1)\Bigl(\frac{n(4d+1)}{n_{\mathrm{out}}}\Bigr)^{n_{\mathrm{out}}/2},
\end{equation}
where the last factor is read as $1$ if $n_{\mathrm{out}}=0$.
\end{lem}
\begin{proof}
The matrix $2A+I$ has diagonal entries $1$ and off-diagonal entries in $\{0,2\}$, so $|\det(2A+I)|\ge1$ by Fact~\ref{fact:parity}. By using the well-known identities for $\operatorname{tr}$ and $\det$ (Fact~\ref{fact:trace}) and~\eqref{eq:input2},
\begin{equation}\label{eq:detidentities}
\det(2A+I)=(2d+1)\prod_{\lambda\in\sigma^*}(2\lambda+1)\qquad\text{and}\qquad\sum_{\lambda\in\sigma^*}(2\lambda+1)^2=4(nd-d^2)-4d+n-1<n(4d+1)\;.
\end{equation}
So the AM-GM inequality (Fact~\ref{fact:amgm}), applied to the numbers $(2\lambda+1)^2$ with $\lambda\in\sigma_{\mathrm{out}}$, gives
\begin{align*}
1\ \le\ |\det(2A+I)|\ &\eqBy{\eqref{eq:detidentities}}\ (2d+1)\Bigl(\prod_{\lambda\in\sigma_{\mathrm{out}}}(2\lambda+1)^2\Bigr)^{1/2}\prod_{\lambda\in\sigma_{\mathrm{mid}}}|2\lambda+1|\\
&\leBy{\ref{fact:amgm}}\ (2d+1)\Bigl(\frac{1}{n_{\mathrm{out}}}\sum_{\lambda\in\sigma_{\mathrm{out}}}(2\lambda+1)^2\Bigr)^{n_{\mathrm{out}}/2}\prod_{\lambda\in\sigma_{\mathrm{mid}}}|2\lambda+1|\\
&\leBy{\eqref{eq:detidentities}}\ (2d+1)\Bigl(\frac{n(4d+1)}{n_{\mathrm{out}}}\Bigr)^{n_{\mathrm{out}}/2}\prod_{\lambda\in\sigma_{\mathrm{mid}}}|2\lambda+1|\;.
\end{align*}
\end{proof}

\subsection*{The errors on the middle cluster are large}
Let
\begin{equation}\label{eq:rho}
\rho\coloneqq\Bigl(\frac{1}{n_{\mathrm{mid}}}\sum_{\lambda\in\sigma_{\mathrm{mid}}}\mu_\lambda^2\Bigr)^{1/2}
\end{equation}
be the root mean square of the errors on the middle cluster ($\rho\coloneqq0$ if $n_{\mathrm{mid}}=0$). Since $|\mu_\lambda|\le r$ by~\eqref{eq:input1}, we have $\rho\le r$. The next proposition makes steps (4) and (5) of the Outline precise, in the worst case allowed by Lemma~\ref{lem:counting}. For real numbers $\eta$ and $n$ with $0\le\eta<1-\frac{2}{n}$ we write
\begin{equation}\label{eq:NM}
N\coloneqq\frac{n}{2}(1+\eta)\qquad\text{and}\qquad M\coloneqq n-1-N=\frac{n}{2}\Bigl(1-\eta-\frac{2}{n}\Bigr)>0
\end{equation}
and
\[
\Phi(d,\eta,n)\coloneqq(8z-4)\Bigl(\frac{1+\eta}{8d+2}\Bigr)^{N/(2M)}(2d+1)^{-1/M}.
\]
So if $n$ is the number of vertices of $G$ and $n_{\mathrm{out}}\le N$, then $M$ is a lower bound for $n_{\mathrm{mid}}=n-1-n_{\mathrm{out}}$.

\begin{prop}[Large Errors]\label{prop:key}
Let $0\le\eta<1-\frac{2}{n}$ be such that $n_{\mathrm{out}}\le N$, for instance as in Lemma~\ref{lem:counting}. Then
\[
1+8\rho\ \ge\ \Phi(d,\eta,n).
\]
Moreover, $\Phi(d,\eta,n)$ is increasing in $n$.
\end{prop}
\begin{proof}
All products and sums in this proof run over $\lambda\in\sigma_{\mathrm{mid}}$, and $n_{\mathrm{mid}}=n-1-n_{\mathrm{out}}\ge n-1-N=M>0$. By the AM-GM inequality (Fact~\ref{fact:amgm}) and the Cauchy--Schwarz inequality $\bigl(\sum|\mu_\lambda|\bigr)^2\le n_{\mathrm{mid}}\sum\mu_\lambda^2$ (CS),
\begin{align}
\prod\bigl(1+8|\mu_\lambda|\bigr)\ &\leBy{\ref{fact:amgm}}\ \Bigl(\frac{1}{n_{\mathrm{mid}}}\sum\bigl(1+8|\mu_\lambda|\bigr)\Bigr)^{n_{\mathrm{mid}}}\ =\ \Bigl(1+8\cdot\frac{1}{n_{\mathrm{mid}}}\sum|\mu_\lambda|\Bigr)^{n_{\mathrm{mid}}}\notag\\
&\leBy{$\mathrm{CS}$}\ \Bigl(1+8\Bigl(\frac{1}{n_{\mathrm{mid}}}\sum\mu_\lambda^2\Bigr)^{1/2}\Bigr)^{n_{\mathrm{mid}}}\ =\ (1+8\rho)^{n_{\mathrm{mid}}}\;.\label{eq:muprod}
\end{align}
As $\eta<1$, the base $\frac{1+\eta}{8d+2}$ in $\Phi$ is less than $1$. Hence, if $1+8\rho\ge8z-4$, the claim holds, as the other two factors of $\Phi$ are at most $1$. Otherwise, by using Lemmas~\ref{lem:clusters} and~\ref{lem:det},
\begin{align*}
\Bigl(\frac{8z-4}{1+8\rho}\Bigr)^{M}\ &\le\ \Bigl(\frac{8z-4}{1+8\rho}\Bigr)^{n_{\mathrm{mid}}}\ \leByRef{eq:muprod}\ \prod\frac{8z-4}{1+8|\mu_\lambda|}\ \leBy{\ref{itm:middlecl}}\ \prod\frac{1}{|2\lambda+1|}\\
&\leByRef{eq:det}\ (2d+1)\Bigl(\frac{n(4d+1)}{n_{\mathrm{out}}}\Bigr)^{n_{\mathrm{out}}/2}\ \le\ (2d+1)\Bigl(\frac{n(4d+1)}{N}\Bigr)^{N/2}\\
&=\ (2d+1)\Bigl(\frac{8d+2}{1+\eta}\Bigr)^{N/2}.
\end{align*}
In the fifth step we used $n_{\mathrm{out}}\le N<n$ and that $x\mapsto\bigl(n(4d+1)/x\bigr)^{x/2}$ is increasing on $(0,n]$. Indeed, the derivative of its logarithm is $\frac{1}{2}\bigl(\log\frac{n(4d+1)}{x}-1\bigr)>0$, as $\frac{n(4d+1)}{x}\ge4d+1>3$. Taking the $M$-th root and rearranging gives the claim. Finally, $\Phi(d,\eta,n)$ is increasing in $n$, as its base $\frac{1+\eta}{8d+2}$ is less than $1$ and its exponents $\frac{N}{2M}=\frac{1+\eta}{2(1-\eta-2/n)}$ and $\frac{1}{M}$ decrease as $n$ increases.
\end{proof}

Since $\frac{N}{2M}\to\frac{1}{2}$ and $\frac{1}{M}\to0$, we have $\Phi\approx\sqrt{8d}$. Explicitly:

\begin{lem}\label{lem:phi}
Let $d\ge100$, $n\ge d^3-d^2$ and $0\le\eta\le1.6/d$. Then $\Phi(d,\eta,n)>\sqrt{8d}-4$.
\end{lem}
\begin{proof}
By~\eqref{eq:NM}, $\frac{N}{2M}+\frac{1}{M}=\frac{1}{2}+\frac{N-M+2}{2M}=\frac{1}{2}+\frac{n\eta+3}{2M}$ and $2M=n(1-\eta)-2$. As $\eta\le1.6/d\le0.016$ and $n\ge d^3-d^2$, we have $2M\ge0.98\,n$ and $\frac{3}{n}\le\frac{0.01}{d}$, so
\begin{equation}\label{eq:phiexp}
\frac{n\eta+3}{2M}\ \le\ \frac{\eta+3/n}{0.98}\ \le\ \frac{1.7}{d}\;.
\end{equation}
Moreover, as $x\mapsto x-\frac{14}{x}$ is increasing,
\begin{equation}\label{eq:phiroot}
\sqrt{8d}-\frac{14}{\sqrt{8d}}\ \le\ \sqrt{8d+2}-\frac{14}{\sqrt{8d+2}}\ =\ \frac{8d-12}{\sqrt{8d+2}}\ \le\ \frac{8d}{\sqrt{8d}}\ =\ \sqrt{8d}\;.
\end{equation}
Hence by using $8z-4=8d-12$, $(1+\eta)^{N/(2M)}\ge1$ and $2d+1<8d+2$,
\begin{align*}
\Phi\ &\ge\ (8d-12)\,(8d+2)^{-\frac{N}{2M}}\,(8d+2)^{-\frac{1}{M}}\ \eqByRef{eq:NM}\ \frac{8d-12}{\sqrt{8d+2}}\,(8d+2)^{-\frac{n\eta+3}{2M}}\\
&\ge\ \frac{8d-12}{\sqrt{8d+2}}\Bigl(1-\frac{n\eta+3}{2M}\log(8d+2)\Bigr)\ \geByRef{eq:phiexp}\ \frac{8d-12}{\sqrt{8d+2}}-\frac{8d-12}{\sqrt{8d+2}}\cdot\frac{1.7\log(8d+2)}{d}\\
&\geByRef{eq:phiroot}\ \sqrt{8d}-\frac{14}{\sqrt{8d}}-\sqrt{8d}\cdot\frac{1.7\log(8d+2)}{d}\;.
\end{align*}
In the third step $y^{-x}=\exp(-x\log y)\ge1-x\log y$ for $x\ge0$ and $y>1$, and in the last step the lower bound in~\eqref{eq:phiroot} for the first term and the upper bound for the second. The two subtracted terms $\frac{14}{\sqrt{8d}}$ and $\frac{1.7\sqrt{8}\,\log(8d+2)}{\sqrt d}$ decrease in $d$, and their sum is less than $3.8$ at $d=100$. Hence $\Phi(d,\eta,n)>\sqrt{8d}-4$.
\end{proof}

\begin{proof}[Proof of Theorem~\ref{thm:deviation7}]
Let $r\coloneqq r(G)$. If $r=0$, then every pair of vertices is joined by exactly one non-backtracking walk of length at most $3$, so $G$ is a Moore graph of diameter $3$. Hence $d\le2$~\cite{HoSi60} and $d<4=(\sqrt8\,r+2)^2$. So let $r\ge1$, as $r$ is an integer, and suppose for a contradiction that $d\ge(\sqrt{8}\,r+2)^2$. Then $d\ge(\sqrt8+2)^2>23$ and $8r\le\sqrt{8}(\sqrt d-2)\le d-3$, so $z\ge8r+2\ge\max(2r+2,3)$. Hence the assumptions of the Setup hold with $r=r(G)$, and moreover $\eps=\frac{8r+1}{8d-12}<\frac{1}{8}$.
By~\eqref{eq:sbound}, $s(G)\le r^2\le\frac{d}{8}$, and by~\eqref{eq:input1}, $n\ge m-r\ge d^3-d^2$. Let $\eta$ be as in Lemma~\ref{lem:counting} with $s=r^2$. By that lemma, $\eta<\frac{1}{10}<1-\frac{2}{n}$.

Let $d\ge100$. Then $\eta\le1.6/d$ by Lemma~\ref{lem:counting}, and
\[
\sqrt{8d}-4\ \lBy{\ref{lem:phi}}\ \Phi(d,\eta,n)\ \leBy{\ref{prop:key}}\ 1+8\rho\ \le\ 1+8r\ \le\ 1+\sqrt{8}\bigl(\sqrt{d}-2\bigr)\ =\ \sqrt{8d}-\bigl(4\sqrt{2}-1\bigr)\;,
\]
a contradiction, as $4\sqrt2-1>4$. For $24\le d\le99$ we have $r\in\{1,2\}$. Then $\Phi(d,\eta,m-r)$ is an explicit expression in $d$ and $r$. By definition,
\[
\log\Phi(d,\eta,n)\ =\ \log(8d-12)-\frac{N}{2M}\log\frac{8d+2}{1+\eta}-\frac{1}{M}\log(2d+1)\;,
\]
where $N$ and $M$ are given by~\eqref{eq:NM} with $n=m-r$, and $\eta=\frac{5/2+\eps+\sqrt{2}\,r/\sqrt{z}}{2d-5/2}$ with $\eps=\frac{8r+1}{8d-12}$. Evaluating this with a computer for the $117$ pairs $(d,r)$ shows $\Phi(d,\eta,m-r)>1+8r$ in each case. As $n\ge m-r$ by~\eqref{eq:input1} and $\Phi$ is increasing in $n$ by Proposition~\ref{prop:key}, we get
\[
1+8\rho\ \le\ 1+8r\ <\ \Phi(d,\eta,m-r)\ \le\ \Phi(d,\eta,n)\,,
\]
contradicting Proposition~\ref{prop:key} as well.
\end{proof}

\subsection*{Half of the error mass is outer}
Theorem~\ref{thm:deviation7} bounds $\rho$ by the largest error $r$. If $C$ or $-C$ is a $0/1$ matrix with row sums $r$, the errors are much smaller on average. In that case $s(G)=r$, so Lemma~\ref{lem:errors}\ref{itm:eigenvands(G)} gives $\sum_{\lambda\in\sigma^*}\mu_\lambda^2=nr-r^2$ instead of the worst case $(n-1)r^2$. As at least about half of the non-trivial eigenvalues are middle, this alone gives $\rho^2\lesssim2r$ and a linear bound $d=O(r)$. The next lemma shows that at least about half of the error mass sits on the outer clusters, so that in fact $\rho^2\lesssim r$; this halves the constant.

\begin{lem}\label{lem:halfmass}
Suppose that $r=r(G)$ and that $C$ or $-C$ has all entries in $\{0,1\}$. Then $|n-m|=r$, $s(G)=r$ and
\[
\sum_{\lambda\in\sigma_{\mathrm{mid}}}\mu_\lambda^2\ \le\ \frac{nr}{2}\Bigl(1+\frac{\tau}{z}\Bigr)+r^2d\;.
\]
\end{lem}
\begin{proof}
Let $D\in\{C,-C\}$ be the $0/1$ matrix.
Recall the polynomial $Q(x)=x^2+x/2$.
By~\eqref{eq:defCk} and Lemma~\ref{lem:polynomials}, every row of $C$ sums to $m-n$. As all entries of $D$ are non-negative, every row of $D$ sums to $|n-m|$ and by~\eqref{eq:defCk} this is also $\sum_v|\nb_3(u,v)-1|$ for every vertex $u$. Therefore $|n-m|=r(G)=r$ and $D\mathbf1=r\mathbf1$. As all entries of $D$ are $0$ or $1$, moreover $s(G)=\frac{1}{n}\operatorname{tr}C^2=\frac{1}{n}\sum_{u,v}D_{u,v}=r$.
Furthermore, $Q(A)D^2$ has all row sums equal to $Q(d)r^2$, so $\mathbf1$ is an eigenvector with eigenvalue $Q(d)r^2$. All other eigenvalues are of the form $Q(\lambda)\mu_\lambda^2$ with $\lambda\in\sigma^*$, since $Q(A)D^2b=\mu_\lambda^2Q(A)b=Q(\lambda)\mu_\lambda^2b$ for the eigenvectors $b\perp\mathbf1$ of Lemma~\ref{lem:errors}\ref{itm:eigenvofC}, and together with $\mathbf1$ these form a basis. Therefore
\begin{equation}\label{eq:trQD2}
\operatorname{tr}\bigl(Q(A)D^2\bigr)=Q(d)\,r^2+\sum_{\lambda\in\sigma^*} Q(\lambda)\mu_\lambda^2\qquad\text{and}\qquad\sum_{\lambda\in\sigma^*}\mu_\lambda^2=\operatorname{tr}D^2-r^2\le\operatorname{tr}D^2=nr\;,
\end{equation}
where the second identity is the same argument with $I$ in place of $Q(A)$ and $\operatorname{tr}(D^2)=nr$, as $D$ is $r$-regular. Since $A$ commutes with $F(A)$ and $J$, it also does with $D$, and hence the matrix $AD$ is symmetric with non-negative integer entries. So $\operatorname{tr}(A^2D^2)=\sum_{u,w}(AD)_{u,w}^2\ge\sum_{u,w}(AD)_{u,w}=\mathbf1^{\mathsf T}AD\mathbf1=d\mathbf1^{\mathsf T}\cdot r\mathbf 1=ndr$, and $\operatorname{tr}(AD^2)\ge0$. This gives
\begin{equation}\label{eq:trQD2low}
\operatorname{tr}(Q(A)D^2)=\operatorname{tr}((A^2+\tfrac{1}{2}A)D^2)\ge ndr\;.
\end{equation}
On the other hand, $Q(\lambda)\le2z+\frac{1}{2}+\tau$ for outer and $Q(\lambda)\le\frac{\eps}{2}\le\frac{1}{2}$ for middle $\lambda$ by Lemma~\ref{lem:clusters}. Writing $X\coloneqq\sum_{\lambda\in\sigma_{\mathrm{mid}}}\mu_\lambda^2$, we obtain
\[
ndr-\Bigl(d^2+\frac{d}{2}\Bigr)r^2\ \By{\eqref{eq:trQD2},\eqref{eq:trQD2low}}{\le}\ \sum_{\lambda\in\sigma^*} Q(\lambda)\mu_\lambda^2\ \By{\eqref{eq:trQD2},\ref{lem:clusters}}{\le}\ \Bigl(2z+\frac{1}{2}+\tau\Bigr)(nr-X)+\frac{X}{2}\;.
\]
Rearranging and using $2z+\frac{1}{2}+\tau-d=z-\frac{1}{2}+\tau$ and $d^2+\frac{d}{2}\le2zd$ gives
\[
2zX\le\ (2z+\tau)X \le\ \Bigl(z-\frac{1}{2}+\tau\Bigr)nr+\Bigl(d^2+\frac{d}{2}\Bigr)r^2\ \le\ (z+\tau)\,nr+2zd\,r^2\;,
\]
which is the claim.
\end{proof}

We now combine Lemma~\ref{lem:halfmass} with Proposition~\ref{prop:key} to prove Theorem~\ref{thm:deviationref7}.

\begin{proof}[Proof of Theorem~\ref{thm:deviationref7}]
Let $r\coloneqq r(G)$. If $r=0$, then $d\le2<4$ as in the proof of Theorem~\ref{thm:deviation7}. So let $r\ge1$ and suppose for a contradiction that $d\ge(\sqrt{8r}+2)^2$. Then $d\ge24$ and $8r\le(\sqrt d-2)^2\le d-3$, so $z\ge8r+2\ge\max(2r+2,3)$, and the assumptions of the Setup hold with $r=r(G)$, and $\eps<\frac{1}{8}$. By Lemma~\ref{lem:halfmass} we have $s(G)=r$ and $|n-m|=r$. So $n\ge m-r\ge d^3-d^2$, as $8r\le d-3$.
Let $\eta$ be as in Lemma~\ref{lem:counting} with $s=r$. By that lemma, $\eta<\frac{1}{10}<1-\frac{2}{n}$ and $n_{\mathrm{out}}\le N$, so $n_{\mathrm{mid}}\ge M>0$. Hence
\begin{equation}\label{eq:rhoref}
\begin{split}
\rho^2\ &\eqByRef{eq:rho}\ \frac{1}{n_{\mathrm{mid}}}\sum_{\lambda\in\sigma_{\mathrm{mid}}}\mu_\lambda^2\ \le\ \frac{1}{M}\sum_{\lambda\in\sigma_{\mathrm{mid}}}\mu_\lambda^2\ \leBy{\ref{lem:halfmass}}\ \frac{1}{M}\Bigl(\frac{nr}{2}\Bigl(1+\frac{\tau}{z}\Bigr)+r^2d\Bigr)\\
&\eqByRef{eq:NM}\ r\cdot\frac{1+\tau/z+2rd/n}{1-\eta-2/n}\;.
\end{split}
\end{equation}
Let $d\ge100$, so that $z\ge0.99\,d$ and $\frac{1}{d}\le\frac{0.1}{\sqrt d}$. As $8r\le d-3<d$, $\tau=\frac{r}{\sqrt{2z}}$ and $n\ge d^3-d^2=d^2z$, we have
\[
\frac{\tau}{z}+\frac{2rd}{n}\ \le\ \frac{d}{8z\sqrt{2z}}+\frac{2d^2}{8d^2z}\ \le\ \frac{1}{8\cdot0.99\sqrt{1.98}\,\sqrt d}+\frac{1}{3.96\,d}\ <\ \frac{0.09+0.026}{\sqrt d}\ <\ \frac{0.12}{\sqrt d}\;.
\]
Moreover, $\eta\le\frac{1.6}{d}$ by Lemma~\ref{lem:counting}, so, using also the crude bound $dz>200$,
\[
\eta+\frac{2}{n}\ \le\ \frac{1.6}{d}+\frac{2}{d^2z}\ \le\ \frac{1.61}{d}\ \le\ \frac{0.17}{\sqrt d}\;.
\]
So, as $\frac{1}{\sqrt d}\le0.1$, the fraction in~\eqref{eq:rhoref} is at most $\bigl(1+\frac{0.12}{\sqrt d}\bigr)/\bigl(1-\frac{0.17}{\sqrt d}\bigr)=\ 1+\frac{0.29}{\sqrt d\,\bigl(1-\frac{0.17}{\sqrt d}\bigr)}\le1+\frac{0.3}{\sqrt d}$, and $\rho\le\sqrt r\,\bigl(1+\frac{0.15}{\sqrt d}\bigr)$.
With $8\sqrt r\le\sqrt8(\sqrt d-2)=\sqrt{8d}-4\sqrt2$ this gives
\[
1+8\rho\ \le\ 1+\bigl(\sqrt{8d}-4\sqrt{2}\bigr)\Bigl(1+\frac{0.15}{\sqrt d}\Bigr)\ \le\ \sqrt{8d}-4\sqrt{2}+1+0.15\sqrt{8}\ <\ \sqrt{8d}-4\ \lBy{\ref{lem:phi}}\ \Phi(d,\eta,n)\;,
\]
contradicting Proposition~\ref{prop:key}. For $24\le d\le99$, the assumption $d\ge(\sqrt{8r}+2)^2$ leaves the $293$ pairs $(d,r)$ with $1\le r\le\frac{(\sqrt d-2)^2}{8}$, and $n=m\pm r$ by Lemma~\ref{lem:halfmass}. So for each of these $293$ pairs $(d,r)$ and both signs, the right-hand side of~\eqref{eq:rhoref} and $\Phi(d,\eta,n)$ are explicit numbers, and evaluating them with a computer shows $1+8\rho<\Phi(d,\eta,n)$ in all cases.
\end{proof}

\section{Girth eleven and diameter five}\label{sec:eleven}

In this section we prove Theorems~\ref{thm:deviation11} and~\ref{thm:deviationref11}. Since only the exponents matter here, we do not track constants.

\subsection*{Setup}
Throughout this section, $G$ is a $d$-regular graph on $n$ vertices with adjacency matrix $A$, and $k=5$. We write $z\coloneqq d-1$ and abbreviate
\[
F\coloneqq F_{d,5}=x^5+x^4-4zx^3-3zx^2+3z^2x+z^2\qquad\text{and}\qquad m\coloneqq m(d,5),
\]
as well as $r(G)\coloneqq r_5(G)$ and $s(G)\coloneqq s_5(G)$. Moreover, $r$ is a real number such that
\begin{equation}\label{eq:rsmall11}
r(G)\ \le\ r\ \le\ \frac{z}{2}\;,
\end{equation}
and $z\ge z_0$ for a sufficiently large absolute constant $z_0$. We write $X=O(Y)$ if $|X|\le KY$ for an absolute constant $K$; these constants do not depend on $z_0$, which is chosen after them. In the proofs of Theorems~\ref{thm:deviation11} and~\ref{thm:deviationref11} we take $r=r(G)$.

As before, $\sigma^*$ is the multiset of the $n-1$ non-trivial eigenvalues of $A$, and $\mu_\lambda=F(\lambda)$ is the error of $\lambda\in\sigma^*$. We also use the rescaled eigenvalues $y_\lambda\coloneqq\lambda/\sqrt z$. By Lemma~\ref{lem:errors} with $k=5$, and as $m\ge dz^4\ge z^5+r$,
\begin{equation}\label{eq:input11}
|\mu_\lambda|\le r\quad\text{for all }\lambda\in\sigma^*,\qquad n\ \ge\ m-r\ \ge\ z^5\qquad\text{and}\qquad\sum_{\lambda\in\sigma^*}\mu_\lambda^2\ \le\ n\,s(G)\;.
\end{equation}

\subsection*{Outline}
The plan is the same as in Section~\ref{sec:seven}.
\begin{enumerate}[label=\textup{(\arabic*)}]
\item By Lemma~\ref{lem:five}, the non-trivial eigenvalues lie in five short clusters, near $-\frac{1}{3}$, $\pm\sqrt z$ and $\pm\sqrt{3z}$.
\item The second and fourth moments of the spectrum determine $\sum_{\lambda\in\sigma^*}q(y_\lambda)$ for every even polynomial $q$ of degree $4$ (Lemma~\ref{lem:spectrum11}). In effect, the five clusters carry about $\frac{n}{12},\frac{n}{4},\frac{n}{3},\frac{n}{4},\frac{n}{12}$ eigenvalues.
\item The middle cluster is pinned near $-\frac{1}{3}$ and the inner clusters near the roots of an integer quadratic $Q$, both up to $O(1/z)$ (Lemma~\ref{lem:pin}). An integer polynomial $P_4$, built from $3x+1$ and $Q$, is therefore small on the middle and the two inner clusters, and large only on the two outer ones (Lemma~\ref{lem:P4}).
\item The integer $\det P_4(A)$ is non-zero (Lemma~\ref{lem:detP4}), so the errors on the middle and inner clusters cannot all be small (Proposition~\ref{prop:key11}).
\end{enumerate}
The only new difficulty is step (3). The linear polynomial $3x+1$, the analogue of $2x+1$, is no longer good enough. We explain this, and how to find a better polynomial, before Lemma~\ref{lem:P4}.

\subsection*{Five clusters}
We shall use two other ways of writing $F$; both are verified by expanding:
\begin{align}
F(x)&=z^2(3x+1)-z\,x^2(4x+3)+x^4(x+1),\label{eq:F5a}\\
F(\sqrt z\,y)&=z^{5/2}f(y)+z^2g(y),\quad\text{where } f(y)\coloneqq y(y^2-1)(y^2-3)\text{ and }g(y)\coloneqq y^4-3y^2+1.\label{eq:F5b}
\end{align}
By~\eqref{eq:F5a}, for bounded $x$ the dominant term of $F(x)$ is $z^2(3x+1)$; this is why one cluster is pinned near $-\frac{1}{3}$. By~\eqref{eq:F5b}, after rescaling by $\sqrt z$ the polynomial $F$ is close to $z^{5/2}f$, whose roots are $0,\pm1,\pm\sqrt3$.

\begin{lem}[Cluster Lemma]\label{lem:five}
Let $|t|\le r$. Every real root $\theta$ of $F(x)=t$ satisfies $|\theta|<2\sqrt z$ and exactly one of the following, where $y\coloneqq\theta/\sqrt z$:
\begin{enumerate}
\labelleditem{(M)}{itm:mid11} $|3\theta+1|=O(1/z)$ \quad (middle cluster, $\theta\approx-\frac{1}{3}$);
\labelleditem{(I)}{itm:in11} $|y^2-1|=O(z^{-1/2})$ \quad (inner clusters, $\theta\approx\pm\sqrt z$);
\labelleditem{(O)}{itm:out11} $|y^2-3|=O(z^{-1/2})$ \quad (outer clusters, $\theta\approx\pm\sqrt{3z}$).
\end{enumerate}
\end{lem}
\begin{proof}
Since $F(\sqrt z\,y)=t$ and $|t|\le r\le z$,~\eqref{eq:F5b} gives
\begin{equation}\label{eq:fsmall}
|f(y)|\ =\ \frac{|t-z^2g(y)|}{z^{5/2}}\ \le\ \frac{|g(y)|+1}{\sqrt z}\;.
\end{equation}
If $|y|\ge2$, then $y^2-1\ge\frac{3}{4}y^2$ and $y^2-3\ge\frac{1}{4}y^2$, so $|f(y)|\ge\frac{3}{16}|y|^5\ge\frac{3}{8}y^4$, whereas $|g(y)|+1\le y^4+1\le2y^4$; this contradicts~\eqref{eq:fsmall} for $z$ large enough. So $|y|<2$, hence $|g(y)|\le5$ and $|f(y)|\le6/\sqrt z$ by~\eqref{eq:fsmall}. Let $y_0\in\{0,\pm1,\pm\sqrt3\}$ be a root of $f$ nearest to $y$. The other four roots have distance at least $\frac{\sqrt3-1}{2}>\frac{1}{3}$ from $y$, half the smallest gap between the roots. As $f(y)=\prod_{\xi}(y-\xi)$, where $\xi$ runs over the five roots of $f$, we get
\[
|y-y_0|=\frac{|f(y)|}{\prod_{\xi\neq y_0}|y-\xi|} \le\ 3^4\,|f(y)|\ \le\ \frac{486}{\sqrt z}\;.
\]
If $y_0=\pm1$, then $|y^2-1|=|y+y_0|\,|y-y_0|\le3|y-y_0|$, which gives~\ref{itm:in11} with the above. If $y_0=\pm\sqrt3$, we get~\ref{itm:out11} in the same way. If $y_0=0$, then $|\theta|=\sqrt z\,|y|\le486$, and~\eqref{eq:F5a} gives
\[
z^2|3\theta+1|\ \le\ |t|+z\,\theta^2|4\theta+3|+\theta^4|\theta+1|\ =\ O(z)\;,
\]
which is~\ref{itm:mid11}. As $z\ge z_0$, the three types are mutually exclusive: they give $y^2\le\frac{1}{z}$, $|y^2-1|\le\frac{1}{2}$ and $|y^2-3|\le\frac{1}{2}$, respectively.
\end{proof}

By~\eqref{eq:input11}, Lemma~\ref{lem:five} applies to every non-trivial eigenvalue $\lambda$, with $\theta=\lambda$ and $t=\mu_\lambda$. We call $\lambda$ \emph{middle}, \emph{inner} or \emph{outer} according to its type, and denote by $\sigma_{\mathrm{mid}}$, $\sigma_{\mathrm{in}}$ and $\sigma_{\mathrm{out}}$ the corresponding sub-multisets of $\sigma^*$, and by $n_{\mathrm{mid}}$, $n_{\mathrm{in}}$ and $n_{\mathrm{out}}$ their sizes. Figure~\ref{fig:five} shows the five clusters.

\begin{figure}[ht]
\centering
\begin{tikzpicture}[x=3.0cm, y=0.5cm]
\foreach \p in {-1.7320508,-1,0,1,1.7320508} \draw[black!22, densely dotted] (\p,4.7) -- (\p,-5.0);
\draw[-{Stealth}, black!60] (-1.95,0) -- (2.05,0) node[right] {$y$};
\draw[thick, smooth, samples=200, domain=-1.86:1.86] plot (\x, {\x*(\x*\x-1)*(\x*\x-3)});
\draw[thick, dashed, smooth, samples=200, domain=-1.86:1.86] plot (\x, {\x*\x*(\x*\x-1)/2});
\node[right, font=\small] at (1.86,2.10) {$f$};
\node[right, font=\small] at (1.86,4.25) {$p$};
\foreach \p in {-1.7320508,-1,0,1,1.7320508} \fill (\p,0) circle (2.3pt);
\foreach \p/\lab in {-1.7320508/{$-\sqrt{3z}$}, -1/{$-\sqrt z$}, 0/{$-\frac{1}{3}$}, 1/{$\sqrt z$}, 1.7320508/{$\sqrt{3z}$}}
  \node[font=\small] at (\p,-5.9) {\lab};
\foreach \p/\lab in {-1.7320508/{$\frac{1}{12}$}, -1/{$\frac{1}{4}$}, 0/{$\frac{1}{3}$}, 1/{$\frac{1}{4}$}, 1.7320508/{$\frac{1}{12}$}}
  \node[font=\small] at (\p,-7.8) {\lab};
\foreach \p/\lab in {-1.7320508/{$O(z^{2})$}, -1/{$O(\frac{r+1}{\sqrt z})$}, 0/{$O(\frac{r+1}{z})$}, 1/{$O(\frac{r+1}{\sqrt z})$}, 1.7320508/{$O(z^{2})$}}
  \node[font=\small] at (\p,-9.8) {\lab};
\node[anchor=east, font=\small, black!70] at (-2.05,-5.9) {cluster near};
\node[anchor=east, font=\small, black!70] at (-2.05,-7.8) {share $\approx$};
\node[anchor=east, font=\small, black!70] at (-2.05,-9.8) {$|P_4(\lambda)|$};
\end{tikzpicture}
\caption{The five clusters, drawn in the rescaled variable $y=x/\sqrt z$. Solid: $f(y)=y(y^2-1)(y^2-3)$; by~\eqref{eq:F5b} and Lemma~\ref{lem:five}, the non-trivial eigenvalues lie in five very short intervals around its roots (dots). Dashed: $p(y)=\frac{1}{2}y^2(y^2-1)$, the leading part of $P_4(\sqrt z\,y)/(18432\,z^2)$ for the auxiliary polynomial $P_4$ of~\eqref{eq:P4}. The polynomial $P_4$ is small on the three inner clusters, which carry about $\frac{5}{6}$ of the spectrum, and large only on the two outer clusters, which carry about $\frac{1}{6}$. (For illustration we give the individual cluster shares; we do not need them in the proof.)}
\label{fig:five}
\end{figure}

\subsection*{Moments}
As for girth $7$, the second moment of the spectrum is $\operatorname{tr}A^2=nd$. Here we also need the fourth moment, and a small deviation $r$ from a Moore graph already excludes $4$-cycles. The two moments determine how many eigenvalues each type of cluster contains, but we only need this for the combined number
\[
N\coloneqq n_{\mathrm{mid}}+n_{\mathrm{in}}=n-1-n_{\mathrm{out}}
\]
of middle and inner eigenvalues.

\begin{lem}[Cluster Size Lemma]\label{lem:spectrum11}
The graph $G$ has no $4$-cycle,
\begin{equation}\label{eq:moments11}
\sum_{\lambda\in\sigma^*}y_\lambda^2\ =\ n+O\Bigl(\frac{n}{z}\Bigr)\qquad\text{and}\qquad\sum_{\lambda\in\sigma^*}y_\lambda^4\ =\ 2n+O\Bigl(\frac{n}{z}\Bigr)\;,
\end{equation}
and
\begin{equation}\label{eq:N11}
\frac{n}{2}\ \le\ N\ =\ \frac{5n}{6}\bigl(1+O(z^{-1/2})\bigr)\;.
\end{equation}
\end{lem}
\begin{proof}
Suppose that $u,a,x,b$ is a $4$-cycle. Every neighbour $v$ of $x$ is the end of two different non-backtracking walks from $u$ of length at most $3$. For $v\notin\{a,b\}$ these are $u,a,x,v$ and $u,b,x,v$, for $v=a$ they are $u,a$ and $u,b,x,a$, and similarly for $v=b$. So the $d$ neighbours $v$ of $x$ all satisfy $|\nb_5(u,v)-1|\ge1$, and
\[
r(G)\ =\ \max_{u'\in V(G)}\sum_{v\in V(G)}\bigl|\nb_5(u',v)-1\bigr|\ \ge\ \sum_{v\in V(G)}\bigl|\nb_5(u,v)-1\bigr|\ \ge\ d\ >\ \frac{z}{2}\;,
\]
contradicting~\eqref{eq:rsmall11}.

So two different vertices have at most one common neighbour. Since $A$ is symmetric, $(A^4)_{u,u}=\sum_w(A^2)_{u,w}(A^2)_{w,u}=\sum_w\bigl((A^2)_{u,w}\bigr)^2$. The term $w=u$ equals $d^2$, as $G$ is $d$-regular. For $w\ne u$ the entry $(A^2)_{u,w}\in\{0,1\}$ is the number of common neighbours of $u$ and $w$, so these terms add up to $\sum_{w\ne u}(A^2)_{u,w}=d^2-d$, the number of walks of length $2$ from $u$ that do not return to $u$. Hence $(A^4)_{u,u}=2d^2-d$ and $\operatorname{tr}A^4=nd(2d-1)$. By Fact~\ref{fact:trace}, $\sum_{\lambda\in\sigma^*}\lambda^2=\operatorname{tr}A^2-d^2$ and $\sum_{\lambda\in\sigma^*}\lambda^4=\operatorname{tr}A^4-d^4$. Using $d=z+1$ and $d(2d-1)=2z^2+3z+1$, we obtain
\[
\sum_{\lambda\in\sigma^*}y_\lambda^2\ =\ \frac{nd-d^2}{z}\ =\ n+\frac{n-d^2}{z}\qquad\text{and}\qquad\sum_{\lambda\in\sigma^*}y_\lambda^4\ =\ \frac{nd(2d-1)-d^4}{z^2}\ =\ 2n+\frac{3n}{z}+\frac{n-d^4}{z^2}\;,
\]
and $d^4\le z^5\le n$ by~\eqref{eq:input11}. This proves~\eqref{eq:moments11}.

For~\eqref{eq:N11} we count the outer eigenvalues. The polynomial $y^2(y^2-1)$ vanishes at $0$ and $\pm1$ and equals $6$ at $\pm\sqrt3$. By Lemma~\ref{lem:five}, $y_\lambda^2\in[0,4)$, and $y_\lambda^2$ lies within $O(z^{-1/2})$ of $0$, $1$ or $3$ for middle, inner or outer $\lambda$, respectively. The polynomial $h(x)\coloneqq x(x-1)$ satisfies $|h'(x)|=|2x-1|\le7$ on $[0,4]$, so $|h(x)-h(x_0)|\le7|x-x_0|$ for $x,x_0\in[0,4]$. As $h(0)=h(1)=0$ and $h(3)=6$, we obtain $y_\lambda^2(y_\lambda^2-1)=O(z^{-1/2})$ for middle and inner $\lambda$, and $y_\lambda^2(y_\lambda^2-1)=6+O(z^{-1/2})$ for outer $\lambda$. Hence
\[
6\,n_{\mathrm{out}}+O\bigl(nz^{-1/2}\bigr)\ =\ \sum_{\lambda\in\sigma^*}y_\lambda^2(y_\lambda^2-1)\ \eqByRef{eq:moments11}\ n+O\Bigl(\frac{n}{z}\Bigr)\;,
\]
so $n_{\mathrm{out}}=\frac{n}{6}\bigl(1+O(z^{-1/2})\bigr)$ and $N=n-1-n_{\mathrm{out}}=\frac{5n}{6}\bigl(1+O(z^{-1/2})\bigr)$. In particular $N\ge\frac{n}{2}$, as $z\ge z_0$.
\end{proof}

In the same way one finds $n_{\mathrm{in}}\approx\frac{n}{2}$ and $n_{\mathrm{mid}}\approx\frac{n}{3}$. We shall not need them, as the proof of Proposition~\ref{prop:key11} uses~\eqref{eq:moments11} directly.

\subsection*{Pinning the middle and the inner clusters}
Lemma~\ref{lem:five} locates the clusters only up to $O(1/z)$ and $O(z^{-1/2})$. For the auxiliary polynomial we need to know more precisely where the middle and the inner clusters sit, and the reason is the same as for girth $7$, where $F_{d,3}(-\frac{1}{2})=\frac{1}{8}$ pinned the middle cluster to $-\frac{1}{2}$. Here, by~\eqref{eq:F5a}, $F(-\frac{1}{3})=-\frac{5z}{27}+\frac{2}{243}$ is small compared with $F'(-\frac{1}{3})\approx3z^2$, so the middle cluster is pinned near $-\frac{1}{3}$. For the inner clusters the role of $-\frac{1}{2}$ is played by the two roots $-\frac{1}{4}\pm\sqrt{z+\frac{3}{16}}$ of
\[
Q(x)\coloneqq8x^2+4x-8z-1\;,
\]
and the role of $F_{d,3}(-\frac{1}{2})=\frac{1}{8}$ by the identity
\begin{equation}\label{eq:FQH}
64\,F(x)\ =\ Q(x)\,H(x)+1-5x\;,\qquad\text{where }H(x)\coloneqq8x^3+4x^2-(24z+1)x-8z+1\;,
\end{equation}
which one checks by expanding. So $F=\frac{1-5x}{64}=O(\sqrt z)$ at the roots of $Q$, whereas $|F'|$ is of order $z^2$ there.

\begin{lem}[Pinning Lemma]\label{lem:pin}
Every non-trivial eigenvalue $\lambda$ satisfies
\begin{align}
3\lambda+1\ &=\ \frac{5}{27z}+O\Bigl(\frac{1+|\mu_\lambda|}{z^{2}}\Bigr)
&&\text{if $\lambda$ is middle,}\label{eq:pinmid}\\
Q(\lambda)\ &=\ -\frac{5}{16z}+O\Bigl(\frac{1+|\mu_\lambda|}{z^{3/2}}\Bigr)
&&\text{if $\lambda$ is inner.}\label{eq:pinin}
\end{align}
\end{lem}
\begin{proof}
Let $\lambda$ be middle and $u\coloneqq3\lambda+1$, so $u=O(1/z)$ by Lemma~\ref{lem:five}. Substituting $\lambda=\frac{u-1}{3}$ gives
\[
\lambda^2(4\lambda+3)\ =\ \frac{(u-1)^2(4u+5)}{27}\ =\ \frac{5-6u}{27}+O(u^2)\qquad\text{and}\qquad\lambda^4(\lambda+1)=O(1)\;,
\]
so by~\eqref{eq:F5a} and $zu^2=O(1/z)$,
\[
\mu_\lambda\ =\ F(\lambda)\ =\ z^2u-\frac{5z}{27}+\frac{6zu}{27}+O(zu^2)+O(1)\ =\ u\Bigl(z^2+\frac{2z}{9}\Bigr)-\frac{5z}{27}+O(1)\;.
\]
Solving for $u$ gives
\[
u\ =\ \frac{5z/27+O(1+|\mu_\lambda|)}{z^2+2z/9}\ =\ \frac{5}{27z}+O\Bigl(\frac{1+|\mu_\lambda|}{z^2}\Bigr)\;,
\]
which is~\eqref{eq:pinmid}.

Let $\lambda$ be inner. By~\eqref{eq:FQH} with $x=\lambda$ we have $Q(\lambda)H(\lambda)=5\lambda+64\mu_\lambda-1$. As $\lambda^2=z+O(\sqrt z)$ by Lemma~\ref{lem:five},
\[
H(\lambda)\ =\ 8\lambda(\lambda^2-3z)+4\lambda^2-\lambda-8z+1\ =\ -16z\lambda+O(z)\;,
\]
and $|\lambda|\ge\frac{1}{2}\sqrt z$, so $|H(\lambda)|\ge8z^{3/2}-O(z)\ge4z^{3/2}$. Hence
\[
Q(\lambda)+\frac{5}{16z}\ =\ \frac{16z\,Q(\lambda)H(\lambda)+5H(\lambda)}{16z\,H(\lambda)}\ =\ \frac{16z\,(64\mu_\lambda-1)+O(z)}{16z\,H(\lambda)}\ =\ O\Bigl(\frac{1+|\mu_\lambda|}{z^{3/2}}\Bigr)\;,
\]
which is~\eqref{eq:pinin}.
\end{proof}

\subsection*{Choosing the auxiliary polynomial}
Let $P$ be a polynomial with integer coefficients such that $\det P(A)\ne0$, and suppose that $|P(\lambda)|$ is about $z^{e_{\mathrm{mid}}}$, $z^{e_{\mathrm{in}}}$ and $z^{e_{\mathrm{out}}}$ on the three types of clusters. Let $q(y)=a+by^2+cy^4$ be the even polynomial with $q(0)=e_{\mathrm{mid}}$, $q(\pm1)=e_{\mathrm{in}}$ and $q(\pm\sqrt3)=e_{\mathrm{out}}$. Then $\log|P(\lambda)|\approx q(y_\lambda)\log z$, and as $|\det P(A)|\ge1$,~\eqref{eq:moments11} gives
\[
0\ \le\ \log|\det P(A)|\ \approx\ \log z\sum_{\lambda\in\sigma^*}q(y_\lambda)\ \By{\eqref{eq:moments11}}{\approx}\ n\,(a+b+2c)\log z\ =\ n\Bigl(\frac{e_{\mathrm{mid}}}{3}+\frac{e_{\mathrm{in}}}{2}+\frac{e_{\mathrm{out}}}{6}\Bigr)\log z\;.
\]
So we obtain a contradiction as soon as the weighted sum in brackets is negative. For $3x+1$ this fails. By Lemma~\ref{lem:pin}, $|3\lambda+1|$ is about $\frac{5}{27z}$ on the middle cluster, and it is about $3\sqrt z$ and $3\sqrt{3z}$ on the inner and outer clusters. A similar calculation as above shows then that the exponents are $-1$, $\frac{1}{2}$, $\frac{1}{2}$ and the weighted sum is $0$. (For girth $7$ the same count for $2x+1$ gives $\frac{1}{2}\cdot\frac{1}{2}-\frac{1}{2}=-\frac{1}{4}$, which is the factor $(6e/\sqrt d)^{n/2}$ in the proof sketch.) To do better, the polynomial must also be small on the inner clusters, which suggests the factor $Q$. The product $(3x+1)Q(x)$ is small on the inner clusters, but no longer on the middle one, so its exponents $0$, $-\frac{1}{2}$, $\frac{3}{2}$ still give the weighted sum $0$. A linear factor and a constant correct this:
\begin{equation}\label{eq:P4}
P_4(x)\coloneqq(3x+1)(384x+371)\,Q(x)+360 ,
\end{equation}
that is, $P_4(x)=9216x^4+16584x^3-(9216z-7804)x^2-(11976z+13)x-(2968z+11)$, has the exponents $-1$, $-\frac{1}{2}$ and $2$ and the weighted sum $-\frac{1}{4}$, as the next lemma shows.

\begin{lem}\label{lem:P4}
Every non-trivial eigenvalue $\lambda$ satisfies
\[
|P_4(\lambda)|\ =\ \begin{cases}
O\bigl((1+|\mu_\lambda|)/z\bigr) & \text{if $\lambda$ is middle},\\
O\bigl((1+|\mu_\lambda|)/\sqrt z\bigr) & \text{if $\lambda$ is inner},\\
O(z^{2}) & \text{if $\lambda$ is outer}.
\end{cases}
\]
\end{lem}
\begin{proof}
If $\lambda$ is middle, then $\lambda=-\frac{1}{3}+O(1/z)$ by Lemma~\ref{lem:five}, so $384\lambda+371=243+O(1/z)$ and $Q(\lambda)=Q(-\frac{1}{3})+O(1/z)=-8z+O(1)$. Hence by~\eqref{eq:pinmid} and~\eqref{eq:P4}
\[
P_4(\lambda)\ =\ \Bigl(\frac{5}{27z}+O\Bigl(\frac{1+|\mu_\lambda|}{z^2}\Bigr)\Bigr)\bigl(243+O(1/z)\bigr)\bigl(-8z+O(1)\bigr)+360\ =\ -360+360+O\Bigl(\frac{1+|\mu_\lambda|}{z}\Bigr)\;.
\]
If $\lambda$ is inner, then $\lambda^2=z+O(\sqrt z)$ by Lemma~\ref{lem:five}, so $(3\lambda+1)(384\lambda+371)=1152\lambda^2+1497\lambda+371=1152z+O(\sqrt z)$. Hence by~\eqref{eq:pinin} and~\eqref{eq:P4}
\[
P_4(\lambda)\ =\ \bigl(1152z+O(\sqrt z)\bigr)\Bigl(-\frac{5}{16z}+O\Bigl(\frac{1+|\mu_\lambda|}{z^{3/2}}\Bigr)\Bigr)+360\ =\ -360+360+O\Bigl(\frac{1+|\mu_\lambda|}{\sqrt z}\Bigr)\;.
\]
If $\lambda$ is outer, then $|\lambda|<2\sqrt z$ by Lemma~\ref{lem:five}, and the expanded form of~\eqref{eq:P4} gives $P_4(\lambda)=O(z^2)$.
\end{proof}

\subsection*{The determinant is non-zero}
For girth $7$ the parity of $\det(2A+I)$ showed that this integer is non-zero. Here parity alone is not enough, as $P_4(A)\equiv A+I\pmod 2$ can have even determinant. Instead, parity shows that an eigenvalue of $A$ which is a root of $P_4$ must be an integer, and $P_4$ has no integer roots.

\begin{lem}[Determinant Lemma]\label{lem:detP4}
For every $d$-regular graph, $\det P_4(A)$ is a non-zero integer; in particular $|\det P_4(A)|\ge1$.
\end{lem}
\begin{proof}
As $P_4$ has integer coefficients, $\det P_4(A)$ is an integer, and by Fact~\ref{fact:trace} it is the product of the numbers $P_4(\theta)$ over the eigenvalues $\theta$ of $A$. Suppose that $P_4(\theta)=0$ for an eigenvalue $\theta$, and let $\psi$ be the minimal polynomial of $\theta$ over $\mathbb Q$. Then $\psi$ is monic and divides the characteristic polynomial of $A$, so $\psi$ has integer coefficients by Gauss's lemma (Fact~\ref{fact:gauss}). As $\psi$ also divides $P_4$ and is monic, division with remainder gives $P_4=\psi\,h$ with $h\in\mathbb Z[x]$. By the expanded form of $P_4$, we have $P_4\equiv x+1\pmod 2$, and since $\psi$ is monic, its reduction modulo $2$ still has degree $\deg\psi$. Hence $\deg\psi=1$, that is, $\theta=a$ is an integer. But then~\eqref{eq:P4} gives $(3a+1)(384a+371)(8a^2+4a-8z-1)=-360$, so $(3a+1)(384a+371)$ divides $360$. This is impossible: the product equals $371$ for $a=0$ and $26$ for $a=-1$, and its absolute value is at least $1985$ for all other integers $a$.
\end{proof}

\subsection*{The errors on the middle and inner clusters are large}
Recall that $N=n_{\mathrm{mid}}+n_{\mathrm{in}}$ is the number of middle and inner eigenvalues, and let
\[
\rho\coloneqq\Bigl(\frac{1}{N}\sum_{\lambda\in\sigma_{\mathrm{mid}}\cup\sigma_{\mathrm{in}}}\mu_\lambda^2\Bigr)^{1/2}
\]
be the root mean square of their errors. Since $|\mu_\lambda|\le r$ by~\eqref{eq:input11}, we have $\rho\le r$. The exponents $-1$, $-\frac{1}{2}$ and $2$ of $z$ in Lemma~\ref{lem:P4} are the values at $y=0$, $\pm1$ and $\pm\sqrt3$ of the even polynomial
\[
q(y)\coloneqq\frac{y^4+y^2-4}{4}\;,
\]
and by~\eqref{eq:moments11}
\begin{equation}\label{eq:qsum}
\sum_{\lambda\in\sigma^*}q(y_\lambda)\ =\ \frac{1}{4}\sum_{\lambda\in\sigma^*}y_\lambda^4+\frac{1}{4}\sum_{\lambda\in\sigma^*}y_\lambda^2-(n-1)\ =\ -\frac{n}{4}+O\Bigl(\frac{n}{z}\Bigr)\;.
\end{equation}
So~\eqref{eq:qsum} is the precise form of the heuristic count before~\eqref{eq:P4}.

\begin{prop}[Large Errors]\label{prop:key11}
We have $\log(1+\rho)\ge\frac{3}{10}\log z-O(1)$.
\end{prop}
\begin{proof}
For $\lambda\in\sigma^*$ let $e_\lambda\in\{-1,-\frac{1}{2},2\}$ be the exponent of $z$ in Lemma~\ref{lem:P4} for the cluster of $\lambda$. Then Lemma~\ref{lem:P4} gives
\[
\log|P_4(\lambda)|\ \le\ e_\lambda\log z+\log(1+|\mu_\lambda|)+O(1)
\]
for middle and inner $\lambda$, and the same without the term $\log(1+|\mu_\lambda|)$ for outer $\lambda$. To sum these bounds over $\lambda$, we replace $e_\lambda$ by $q(y_\lambda)$, whose sum we know from~\eqref{eq:qsum}. As a function of $y^2$, we have that the derivative of $q$ is bounded on $[0,4]$. By Lemma~\ref{lem:five}, $y_\lambda^2\in[0,4)$ lies within $O(z^{-1/2})$ of the point $0$, $1$ or $3$ at which $q$ takes the value $e_\lambda$, so $|e_\lambda-q(y_\lambda)|=O(z^{-1/2})$ and
\[
e_\lambda\log z\ \le\ q(y_\lambda)\log z+O\bigl(z^{-1/2}\log z\bigr)\ =\ q(y_\lambda)\log z+O(1)\;.
\]
Hence
\begin{equation}\label{eq:logP4}
\log|P_4(\lambda)|\ \le\ q(y_\lambda)\log z+\log(1+|\mu_\lambda|)+O(1)
\end{equation}
for middle and inner $\lambda$, and the same without the term $\log(1+|\mu_\lambda|)$ for outer $\lambda$. By the AM-GM and the Cauchy--Schwarz inequality, exactly as in~\eqref{eq:muprod},
\begin{equation}\label{eq:muprod11}
\sum_{\lambda\in\sigma_{\mathrm{mid}}\cup\sigma_{\mathrm{in}}}\log(1+|\mu_\lambda|)\ \leBy{\ref{fact:amgm}}\ N\log\Bigl(1+\frac{1}{N}\sum_{\lambda\in\sigma_{\mathrm{mid}}\cup\sigma_{\mathrm{in}}}|\mu_\lambda|\Bigr)\ \leBy{CS}\ N\log(1+\rho)\;.
\end{equation}
Finally, $P_4(d)=O(z^4)$. So by using Lemma~\ref{lem:detP4} and Fact~\ref{fact:trace}
\begin{align*}
0\ \leBy{\ref{lem:detP4}}\ \log|\det P_4(A)|\ &\eqBy{\ref{fact:trace}}\ \log|P_4(d)|+\sum_{\lambda\in\sigma^*}\log|P_4(\lambda)|\\
&\leByRef{eq:logP4}\ O(\log z)+\sum_{\lambda\in\sigma^*}\bigl(q(y_\lambda)\log z+O(1)\bigr)+\sum_{\lambda\in\sigma_{\mathrm{mid}}\cup\sigma_{\mathrm{in}}}\log(1+|\mu_\lambda|)\\
&\leByRef{eq:muprod11}\ \log z\sum_{\lambda\in\sigma^*}q(y_\lambda)+N\log(1+\rho)+O(n)\\
&\leByRef{eq:qsum}\ -\frac{n}{4}\log z+N\log(1+\rho)+O(n)\;.
\end{align*}
So $N\log(1+\rho)\ge n\bigl(\frac{1}{4}\log z-O(1)\bigr)$. As the right-hand side is positive and $N\le\frac{5n}{6}\bigl(1+O(z^{-1/2})\bigr)$ by~\eqref{eq:N11}, we get
\[
\log(1+\rho)\ \geByRef{eq:N11}\ \frac{6}{5}\cdot\frac{\frac{1}{4}\log z-O(1)}{1+O(z^{-1/2})}\ \ge\ \frac{6}{5}\bigl(1-O(z^{-1/2})\bigr)\Bigl(\frac{1}{4}\log z-O(1)\Bigr)\ =\ \frac{3}{10}\log z-O(1)\;,
\]
where we used $\frac{1}{1+x}\ge1-x$ for $x\ge0$ and $z^{-1/2}\log z=O(1)$.
\end{proof}

\begin{proof}[Proof of Theorem~\ref{thm:deviation11}]
Let $K$ be an absolute constant such that Proposition~\ref{prop:key11} holds in the form $\log(1+\rho)\ge\frac{3}{10}\log z-K$, and let $C\coloneqq\max\bigl(z_0+2,3\exp(10K/3)\bigr)$. Let $r=r(G)$. If $z<z_0$, then $d<C$, and if $r>\frac{z}{2}$, then $d<2(r+1)<C(r+1)$. Otherwise~\eqref{eq:rsmall11} holds, and Proposition~\ref{prop:key11} and $\rho\le r$ give $\log(1+r)\ge\frac{3}{10}\log z-K$, that is, $d\le2z\le2\exp(10K/3)\,(r+1)^{10/3}$. In all cases $d<C(r+1)^{10/3}$.
\end{proof}

\begin{proof}[Proof of Theorem~\ref{thm:deviationref11}]
Let $K$ be as in the previous proof, and let $C\coloneqq\max\bigl(z_0+2,21\exp(10K/3)\bigr)$. Let $r=r(G)$ and $s=s(G)$, and suppose that $d\ge C(r+1)$. Then $z=d-1\ge C-1>z_0$ and $r\le\frac{d}{C}-1\le\frac{d}{2}-1<\frac{z}{2}$, so~\eqref{eq:rsmall11} holds, and
\[
\rho^2\ \le\ \frac{1}{N}\sum_{\lambda\in\sigma^*}\mu_\lambda^2\ \By{\eqref{eq:input11},\,\eqref{eq:N11}}{\le}\ \frac{2}{n}\cdot n\,s\ =\ 2s\;.
\]
So Proposition~\ref{prop:key11} gives $\log(1+\sqrt{2s})\ge\frac{3}{10}\log z-K$. As $(1+\sqrt{2s})^2\le4(1+s)$, we get $d\le2z\le2\cdot4^{5/3}\exp(10K/3)\,(s+1)^{5/3}<C(s+1)^{5/3}$.
\end{proof}

\section{Concluding remarks}\label{sec:remarks}

\subsection*{Limits of the method}
Our argument rests on the identity $F_{d,3}(-\frac{1}{2})=\frac{1}{8}$, which holds for every $d$ and pins the middle cluster to the rational number $-\frac{1}{2}$. Considering $2A+I$ shifts this cluster close to zero, and this is exactly what lets $\det(2A+I)$ tend to zero for the final contradiction. For girth $11$ the middle cluster is pinned near $-\frac{1}{3}$, but it carries only a third of the spectrum, so $3A+I$ is not enough and we used the quartic $P_4$ instead; see the discussion before Lemma~\ref{lem:P4}. We did not find a more suitable polynomial in our computer search up to degree $8$, but there are polynomials of degree at least $9$ which could potentially lead to better exponents for the excess and defect in Theorems~\ref{thm:girth11} and~\ref{thm:diam5}. We did not pursue those in favour of the much more convenient polynomial $P_4$.

For $k=2$, that is, for girth $5$ and diameter $2$, no polynomial suitable for our method exists. The roots of $F_{d,2}(x)=x^2+x-(d-1)$ are $\lambda_\pm=\frac{1}{2}(-1\pm\sqrt{4d-3})$, and $\operatorname{tr}A=0$ implies that the two corresponding clusters each carry about half of the spectrum. For the final contradiction, a polynomial $P$ needs
\[
|\det P(A)|\ =\ |P(d)|\prod_{\lambda\in\sigma^*}|P(\lambda)|\ \approx\ |P(d)|\,|P(\lambda_+)|^{n/2}\,|P(\lambda_-)|^{n/2}\ \to\ 0
\]
when all errors are small. So if $|P(\lambda_\pm)|\approx d^{e_\pm}$, we need $e_++e_-<0$. We search for $P$ whose coefficients are polynomials in $d$ with integer coefficients, as these can follow clusters that move with $d$. Let $\alpha +\beta x$ be the remainder of $P$ modulo $F_{d,2}$, so $\alpha$ and $\beta$ are again polynomials in $d$ with integer coefficients. As $\lambda_\pm$ are the roots of $F_{d,2}$, we have $\lambda_++\lambda_-=-1$, $\lambda_+\lambda_-=-(d-1)$ and
\[
P(\lambda_+)P(\lambda_-)\ =\ (\alpha+\beta\lambda_+)(\alpha+\beta\lambda_-)\ =\ \alpha^2-\alpha\beta-(d-1)\beta^2\ =:\ R(d),
\]
a polynomial in $d$ with integer coefficients. If $R\ne0$, then $d^{e_++e_-}\approx|R(d)|$ gives $e_++e_-=\deg R\ge0$, so the determinant does not tend to $0$. If $R=0$, then $\alpha=\beta=0$, so $F_{d,2}$ divides $P$ and $P(\lambda)=\mu_\lambda H(\lambda)$ for the eigenvalues, where $H=P/F_{d,2}$. Then $\det P(A)$ contains the product of the errors, and no integrality argument can bound it away from $0$.

For $k=4$ no such polynomial exists either, although the four clusters carry different shares. The roots of $F_{d,4}$ are series in powers of $\sqrt d$ with coefficients in $\mathbb{Q}(\sqrt5)$ and leading terms $\pm\varphi\sqrt d$ and $\pm\varphi^{-1}\sqrt d$, where $\varphi=\frac{1+\sqrt5}{2}$, and they arise from each other by changing the signs of $\sqrt d$ and of $\sqrt5$. As the coefficients of $P$ are again polynomials in $d$ with integer coefficients, its value at one root is again such a series, and its values at the other three roots are obtained by the same sign changes. In particular, if $P$ vanishes at one root, then it vanishes at all four, so $F_{d,4}$ divides $P$, and we argue as for $k=2$. Otherwise, the sign changes turn the leading coefficient $a+b\sqrt5$ of the value at one root, with $a,b$ rational, into $\pm(a\pm b\sqrt5)$, which is again non-zero as $\sqrt5$ is irrational. Hence all four values have the same order $d^e$. As for $k=2$, their product is a polynomial in $d$ with integer coefficients, so $e\ge0$, and the weighted average of $\log|P|$ over the clusters is about $e\log d\ge0$, whatever the shares. We expect a similar obstruction whenever $k+1$ is prime. If $k+1$ is composite, the clusters split into several orbits with different shares, and more complicated polynomials of higher degree might make a similar argument work. However, this would likely involve lengthy calculations for each $k$, and we do not see how to turn it into a general argument for all large $k$.

\subsection*{Tightness of the proof}
Within this framework our bounds for girth $7$ and diameter $3$ are essentially best possible. The method shows $1+8\rho\gtrsim\sqrt{8d}$ for the root mean square $\rho$ of the errors on the middle cluster. Nothing we use prevents $\rho$ from being about $r_3(G)$ in general, or about $\sqrt e$ for excess $e$, as it would be if the excess graph looked like a random $e$-regular graph. Going beyond would need information on how the errors, in the excess case the eigenvalues of the excess graph, are distributed over the three clusters.

For a general graph $G$ of maximum degree $d$, diameter $3$ and defect $f$, Theorem~\ref{thm:deviationref7} does not apply, as the defect graph (the multigraph with adjacency matrix $C_3$) may have multiple edges. If their multiplicities were bounded, or more generally $s_3(G)=O(f)$, the same argument as for the excess would still give a linear bound for the defect, but we see no way to guarantee this without further structural information. Thus an unconditional linear bound for the defect remains open.

\subsection*{Cages}
A \emph{$(d,g)$-cage} is a $d$-regular graph of girth $g$ with the smallest possible number of vertices; see the survey of Exoo and Jajcay~\cite{ExJa13}. Since a $d$-regular graph has minimum degree $d$, Theorems~\ref{thm:girth7} and~\ref{thm:girth11} apply to cages. The excess of a $(d,7)$-cage is larger than $(\sqrt d-2)^2/8$, and that of a $(d,11)$-cage is larger than $c\,d^{3/5}$. In particular, this answers in the negative the question BB1 of Eze and Jajcay~\cite{EzJa24}, which is the question of Bollob\'as for regular graphs: is there a constant $c$ such that for all $d$ and $g$ there is a $d$-regular graph of girth $g$ whose order exceeds the Moore bound by at most $c$. Our method fixes the girth and lets the degree grow. It says nothing about a fixed degree and growing girth, which is the regime of the second part of Problem~11 and of Problem~12 in~\cite{Bol78}, and of the questions BB2 and BB3 in~\cite{EzJa24}.

\subsection*{Comparison with the upper bound}
Known constructions can achieve the following.
For girth $7$, we have $n(d,3)\le2dq^2$ for the smallest odd prime power $q\ge d$~\cite{LUW97}, so that $n(d,3)\le(2+o(1))\,d^3$ by the known bounds on gaps between primes~\cite{BHP01}; for $d\in\{q,q+1\}$ slightly smaller constructions are known~\cite{AAPBLS15,Kis23}. Therefore, together with Theorem~\ref{thm:girth7}, the excess of $(d,7)$-cages lies between $(1-o(1))\,d/8$ and $(1+o(1))\,d^3$, and it remains open whether $n(d,3)/m(d,3)$ tends to $1$. For girth $11$, if $d-1$ is a prime power, the incidence graph of a generalised hexagon of order $d-1$ is $d$-regular, has girth $12$ and $m(d,5)+(d-1)^5$ vertices~\cite{Ben66}, so the excess of $(d,11)$-cages lies between $c\,d^{3/5}$ and $(d-1)^5$ for these degrees.

For the degree/diameter problem the situation is different: $\tilde n(d,k)=(1-o(1))\,d^k$ for every fixed $k$~\cite{CvBK26}. So by Theorems~\ref{thm:diam3} and~\ref{thm:diam5} the defect of extremal graphs of diameter $3$ lies between $\Omega(\sqrt d)$ and $o(d^3)$, and that of diameter $5$ between $\Omega(d^{3/10})$ and $o(d^5)$.

\subsection*{Even girth}
For a graph $G$ of even girth $2k$ and minimum degree $d$, the Moore bound is $\widetilde{m}(d,k)\coloneqq m(d,k-1)+(d-1)^{k-1}$, as a breadth-first search from the two ends of an edge finds at least this many vertices. This bound is attained by the incidence graphs of generalised $k$-gons of order $d-1$. For $k=2$ these are the complete bipartite graphs $K_{d,d}$, and for $k\in\{3,4,6\}$ they exist whenever $d-1$ is a prime power~\cite{Ben66}. For other $d$ their existence is open in general, which for $k=3$ is the famous question about projective planes of non-prime-power order. By a theorem of Feit and Higman~\cite{FeHi64}, the bound cannot be attained for any other $k$ when $d\ge3$; see also~\cite{ExJa13}.

We do not expect our method to apply to even girth, and the reason can be seen at the level of the polynomials. Let $G$ be $d$-regular of girth at least $2k$. Then every non-backtracking walk of length at most $k$ is a path, and no vertex at distance at most $k-1$ from $u$ is the end of a non-backtracking walk of length $k$ from $u$. Hence the matrix
\[
H_{d,k}(A)\coloneqq F_{d,k}(A)+(d-1)F_{d,k-1}(A)
\]
has the entry $d$ at pairs at distance at most $k-1$, the number of paths of length $k$ at pairs at distance $k$, and $0$ otherwise. Its row sums are $d\cdot\widetilde{m}(d,k)$, and for the incidence graph of a generalised $k$-gon we have $H_{d,k}(A)=dJ$, as opposite elements are joined by exactly $d$ paths. So $H_{d,k}$ plays the role of $F_{d,k}$: for the generalised polygons, every eigenvalue of $A$ other than $d$ is a root of $H_{d,k}$. One checks that $-d$ is a root of $H_{d,k}$ for every $k$, as it must be for bipartite graphs, and that (with $z=d-1$)
\begin{align*}
H_{d,2}(x)&=x(x+d), & H_{d,3}(x)&=(x+d)(x^2-z),\\
H_{d,4}(x)&=x(x+d)(x^2-2z), & H_{d,6}(x)&=x(x+d)(x^2-z)(x^2-3z).
\end{align*}
More generally, $0$ is a root of $H_{d,k}$ for every even $k$. So for girths $8$ and $12$ the analogue of the middle cluster is pinned to the integer $0$ itself, not to a rational number with a non-trivial denominator. The integer playing the role of $\det(2A+I)$ would be $\det A$, but nothing forces it to be non-zero. Indeed, the incidence graphs of generalised quadrangles and hexagons have the eigenvalue $0$. For odd $k$, there is not even such a candidate. The roots of $H_{d,k}$ other than $-d$ are $2\sqrt{z}\cos\frac{j\pi}{k}$ with $1\le j\le k-1$, all non-zero multiples of $\sqrt z$, and an obstruction similar to the case of odd girth and $k=4$ can also apply here; for girth $10$, for instance, they are $\pm\varphi\sqrt z$ and $\pm\varphi^{-1}\sqrt z$. Finally, for girths $6$, $8$ and $12$, a proof that the excess is unbounded for $d-1$ not a prime power would in particular exclude projective planes, generalised quadrangles and generalised hexagons.

\section*{AI declaration}
Opus 5.5 was used to assist with the research. Specifically, it was used to discuss the sharpening of some of the bounds in Sections~\ref{sec:seven} and~\ref{sec:eleven} and to write code to search for a suitable auxiliary polynomial for the case of $k=5$. It was also used to improve the presentation of the paper, especially for finding neat expressions to simplify calculations, and to assist with the creation of figures. The author takes full accountability for the text and mathematical content of this paper.

\appendix
\section{Facts from linear algebra}\label{app:facts}

We use the following standard facts, which we record here for completeness.

\begin{fact}\label{fact:spectral}
A real symmetric matrix $M$ has real eigenvalues and an orthonormal basis of eigenvectors. If $b$ is an eigenvector of $M$, then $M$ maps $b^\perp$ into itself, so the basis can be chosen to consist of $b/\norm{b}$ and vectors orthogonal to $b$. Commuting real symmetric matrices have a common orthonormal basis of eigenvectors.
\end{fact}

\begin{fact}\label{fact:rowsum}
If every row of a real matrix $M$ satisfies $\sum_j|M_{i,j}|\le\varrho$, then every real eigenvalue $\mu$ of $M$ satisfies $|\mu|\le\varrho$.
\end{fact}
\begin{proof}
Let $Mb=\mu b$ and choose $i$ with $|b_i|$ maximal. Then $|\mu|\,|b_i|=|\sum_jM_{i,j}b_j|\le\varrho\,|b_i|$.
\end{proof}

\begin{fact}\label{fact:trace}
Let $M$ be a real symmetric matrix with eigenvalues $\theta_1,\dots,\theta_n$. For every polynomial $P$, the matrix $P(M)$ has the eigenvalues $P(\theta_1),\dots,P(\theta_n)$, with the same eigenvectors, since $Mb=\theta b$ implies $P(M)b=P(\theta)b$. Moreover, $\operatorname{tr}M=\sum_i\theta_i$ and $\det M=\prod_i\theta_i$. If $A$ is the adjacency matrix of a $d$-regular graph on $n$ vertices, then $\operatorname{tr}A=0$ and $\operatorname{tr}A^2=nd$, because $(A^2)_{v,v}$ counts the walks $v,w,v$.
\end{fact}

\begin{fact}\label{fact:parity}
If $M$ is an integer matrix whose diagonal entries are odd and whose off-diagonal entries are even, then $\det M$ is an odd integer; in particular $|\det M|\ge1$.
\end{fact}
\begin{proof}
In the Leibniz formula $\det M=\sum_\sigma\operatorname{sgn}(\sigma)\prod_iM_{i,\sigma(i)}$ every term is an integer. Every permutation $\sigma\ne\mathrm{id}$ uses an off-diagonal entry, so its term is even, while the term of $\sigma=\mathrm{id}$ is a product of odd numbers.
\end{proof}

\begin{fact}[AM--GM]\label{fact:amgm}
For non-negative reals $a_1,\dots,a_p$ we have $a_1\cdots a_p\le\bigl(\tfrac1p\sum_ia_i\bigr)^p$.
\end{fact}

\begin{fact}\label{fact:gauss}
\leavevmode
\begin{enumerate}[label=\textup{(\alph*)}]
\item The characteristic polynomial $\det(xI-M)$ of an integer matrix $M$ is monic with integer coefficients, and its roots are the eigenvalues of $M$.
\item If a monic polynomial with rational coefficients divides a monic polynomial with integer coefficients, then it has integer coefficients. This is a standard consequence of Gauss's lemma, which says that the product of two primitive polynomials with integer coefficients is primitive.
\end{enumerate}
\end{fact}

\bibliographystyle{plain}
\bibliography{Moore7_11}

\begin{thebibliography}{10}

\bibitem{AAPBLS15}
Marien Abreu, Gabriela Araujo-Pardo, Camino Balbuena, Domenico Labbate, and
  Juli{\'a}n Salas.
\newblock Small regular graphs of girth 7.
\newblock {\em Electronic Journal of Combinatorics}, 22(3):Paper P3.5, 2015.

\bibitem{ABLS07}
Noga Alon, Itai Benjamini, Eyal Lubetzky, and Sasha Sodin.
\newblock Non-backtracking random walks mix faster.
\newblock {\em Communications in Contemporary Mathematics}, 9(4):585--603,
  2007.

\bibitem{AlHoLi02}
Noga Alon, Shlomo Hoory, and Nathan Linial.
\newblock The {M}oore bound for irregular graphs.
\newblock {\em Graphs and Combinatorics}, 18:53--57, 2002.

\bibitem{BSS19}
Martin Bachrat{\'y}, Jana {\v{S}}iagiov{\'a}, and Jozef {\v{S}}ir{\'a}{\v{n}}.
\newblock Asymptotically approaching the {M}oore bound for diameter three by
  {C}ayley graphs.
\newblock {\em Journal of Combinatorial Theory, Series B}, 134:203--217, 2019.

\bibitem{BHP01}
R.~C. Baker, G.~Harman, and J.~Pintz.
\newblock The difference between consecutive primes, {II}.
\newblock {\em Proceedings of the London Mathematical Society}, 83(3):532--562,
  2001.

\bibitem{BaIt73}
Eiichi Bannai and Tatsuro Ito.
\newblock On finite {M}oore graphs.
\newblock {\em Journal of the Faculty of Science, University of Tokyo, Section
  IA, Mathematics}, 20:191--208, 1973.

\bibitem{BaIt81}
Eiichi Bannai and Tatsuro Ito.
\newblock Regular graphs with excess one.
\newblock {\em Discrete Mathematics}, 37:147--158, 1981.

\bibitem{Ben66}
Clark~T. Benson.
\newblock Minimal regular graphs of girths eight and twelve.
\newblock {\em Canadian Journal of Mathematics}, 18:1091--1094, 1966.

\bibitem{BeBo81}
Jean-Claude Bermond and B{\'e}la Bollob{\'a}s.
\newblock The diameter of graphs---a survey.
\newblock {\em Congressus Numerantium}, 32:3--27, 1981.

\bibitem{Big82}
Norman~L. Biggs.
\newblock Excess in vertex-transitive graphs.
\newblock {\em Bulletin of the London Mathematical Society}, 14:52--54, 1982.

\bibitem{Bol78}
B{\'e}la Bollob{\'a}s.
\newblock {\em Extremal Graph Theory}, volume~11 of {\em London Mathematical
  Society Monographs}.
\newblock Academic Press, London, 1978.

\bibitem{Bro67}
William~G. Brown.
\newblock On the non-existence of a type of regular graphs of girth 5.
\newblock {\em Canadian Journal of Mathematics}, 19:644--648, 1967.

\bibitem{CvBK26}
Wouter {Cames van Batenburg} and Samuel Korsky.
\newblock Asymptotically attaining the {M}oore bound, 2026.
\newblock arXiv:2608.03965.

\bibitem{Dam73}
R.~M. Damerell.
\newblock On {M}oore graphs.
\newblock {\em Proceedings of the Cambridge Philosophical Society},
  74:227--236, 1973.

\bibitem{Del85a}
Charles Delorme.
\newblock Grands graphes de degr{\'e} et diam{\`e}tre donn{\'e}s.
\newblock {\em European Journal of Combinatorics}, 6:291--302, 1985.

\bibitem{Del85b}
Charles Delorme.
\newblock Large bipartite graphs with given degree and diameter.
\newblock {\em Journal of Graph Theory}, 9:325--334, 1985.

\bibitem{ExJa13}
Geoffrey Exoo and Robert Jajcay.
\newblock Dynamic cage survey.
\newblock {\em Electronic Journal of Combinatorics}, 2013.
\newblock Dynamic Survey DS16.

\bibitem{EzJa24}
Leonard~Chidiebere Eze and Robert Jajcay.
\newblock Analogues of {B}ermond--{B}ollob{\'a}s conjecture for cages yield
  expander families.
\newblock {\em Electronic Journal of Graph Theory and Applications},
  14(1):139--152, 2026.

\bibitem{FeHi64}
Walter Feit and Graham Higman.
\newblock The nonexistence of certain generalized polygons.
\newblock {\em Journal of Algebra}, 1(2):114--131, 1964.

\bibitem{FeMiPi11}
Ramiro Feria-Pur{\'o}n, Mirka Miller, and Guillermo Pineda-Villavicencio.
\newblock On graphs of defect at most 2.
\newblock {\em Discrete Applied Mathematics}, 159(13):1331--1344, 2011.

\bibitem{FiJa21}
Slobodan Filipovski and Robert Jajcay.
\newblock A connection between a question of {B}ermond and {B}ollob{\'a}s and
  {R}amanujan graphs.
\newblock {\em Acta Applicandae Mathematicae}, 175:Paper No. 1, 2021.

\bibitem{Gar15}
Frederik Garbe.
\newblock On graphs with excess or defect 2.
\newblock {\em Discrete Applied Mathematics}, 180:81--88, 2015.

\bibitem{HoSi60}
Alan~J. Hoffman and Robert~R. Singleton.
\newblock On {M}oore graphs with diameters 2 and 3.
\newblock {\em IBM Journal of Research and Development}, 4:497--504, 1960.

\bibitem{Kis23}
Gy{\"o}rgy Kiss.
\newblock Geometric constructions of small regular graphs with girth 7.
\newblock {\em The Art of Discrete and Applied Mathematics}, 8(3):Paper P3.04,
  2025.

\bibitem{Kov81}
Peter Kov{\'a}cs.
\newblock The non-existence of certain regular graphs of girth 5.
\newblock {\em Journal of Combinatorial Theory, Series B}, 30(3):282--284,
  1981.

\bibitem{LUW97}
Felix Lazebnik, Vasiliy~A. Ustimenko, and Andrew~J. Woldar.
\newblock New upper bounds on the order of cages.
\newblock {\em Electronic Journal of Combinatorics}, 4(2):Paper R13, 1997.

\bibitem{MiSi13}
Mirka Miller and Jozef {\v{S}}ir{\'a}{\v{n}}.
\newblock Moore graphs and beyond: {A} survey of the degree/diameter problem.
\newblock {\em Electronic Journal of Combinatorics}, 2013.
\newblock Dynamic Survey DS14.

\bibitem{Sin66}
Robert~R. Singleton.
\newblock On minimal graphs of maximum even girth.
\newblock {\em Journal of Combinatorial Theory}, 1(3):306--332, 1966.

\bibitem{VMa98}
Hendrik Van~Maldeghem.
\newblock {\em Generalized Polygons}, volume~93 of {\em Monographs in
  Mathematics}.
\newblock Birkh{\"a}user, Basel, 1998.

\end{thebibliography}

\end{document}